\documentclass[12pt]{article}
\usepackage[utf8]{inputenc}
\usepackage[T1]{fontenc}
\usepackage[left=0.75in,right=0.75in,top=0.95in,bottom=0.95in]{geometry}
\usepackage{mathptmx}
\usepackage{microtype}
\usepackage[style=alphabetic]{biblatex}
\usepackage{amsmath, amssymb, amsthm, stmaryrd, bbm, mathrsfs, url, enumitem, graphicx, xcolor}
\usepackage{times,bm, booktabs}
\usepackage{hyperref}

\hypersetup{
    colorlinks=true,
    citecolor=green!50!black,
    linkcolor=olive
    }

\newtheorem{theorem}{Theorem}
\newtheorem{proposition}{Proposition}
\newtheorem*{lemma*}{Lemma}
\theoremstyle{definition}
\newtheorem{example}{Example}
\theoremstyle{remark}
\newtheorem{remark}{Remark}

\def\thmref#1{\textsc{Theorem}~\ref{#1}}
\def\propref#1{\textsc{Proposition}~\ref{#1}}

\def\secref#1{Section~\ref{#1}}

\def\bfmath#1{\mathchoice
        {\mbox{\boldmath$#1$}}%
        {\mbox{\boldmath$#1$}}%
        {\mbox{\boldmath$\scriptstyle#1$}}%
        {\mbox{\boldmath$\scriptscriptstyle#1$}}}%

\def\discreteS{D}

\def\bfS{\bfmath{S}}

\def\bfe{\bfmath{e}}
\def\bfh{\bfmath{h}}
\def\bfm{\bfmath{m}}
\def\bfr{\bfmath{r}}
\def\bfs{\bfmath{s}}
\def\bfw{\bfmath{w}}
\def\bfx{\bfmath{x}}
\def\iid{\text{i.i.d.}}

\def\sample{\sigma}

\def\lineseg#1#2{\overleftrightarrow{\hspace*{1mm}#1,#2\hspace*{1mm}}}

\begin{document}

\title{Exact and arbitrarily accurate non-parametric two-sample tests based on rank
spacings}
\date{}

\author{
\textsc{Dan D. Erdmann-Pham} \\
\textit{\normalsize Department of Statistics, Stanford University, CA 94305, U.S.A.} \\
{\normalsize erdpham@stanford.edu}
\\[5mm]
\textsc{Jonathan Terhorst} \\
\textit{\normalsize Department of Statistics, University of Michigan, Ann Arbor, MI 48109, U.S.A.} \\
{\normalsize jonth@umich.edu}
\\[5mm]
\textsc{Yun S. Song} \\
\textit{\normalsize Department of Statistics and Computer Science Division, UC Berkeley, CA 94720, U.S.A.} \\
{\normalsize yss@berkeley.edu}
}

\maketitle

\begin{abstract}
    A common method for deriving non-parametric tests is to 
    reformulate a parametric test in terms of sample ranks.
    Despite being distribution free (even
    in finite samples), the resulting tests often display remarkable
    asymptotic power properties, typically matching the efficiency of their
    parametric counterpart. Empirically, these favorable power properties
    have been shown to persist in non-asymptotic regimes as well, prompting
    the need for finite-sample characterizations of the corresponding
    rank-based statistics. Here, we provide such characterization for the
    family of weighted $p$-norms of rank spacings, which includes the
    classical tests of Mann-Whitney, Dixon, and various 
    generalizations thereof. For $p=1$, we provide exact expressions for the
    involved distributions, while for $p>1$ we describe the associated
    moment sequences and derive an algorithm to recover the distributions
    of interest from these sequences in a fast and stable manner. We
    use this framework to develop a new family of non-parametric tests
    mirroring properties of generalized likelihood-ratios, prove new tail
    bounds for Dixon's and Greenwood's statistics, and prove a
    previously formulated conjecture regarding the global efficiency of
    rank-based tests against the $F$-test in the context of scale-families.
\end{abstract}
\newpage
\section{Introduction}
\label{sec:introduction}
Given a pair of samples $\mathscr X_k = \left\{ X_j \right\}_{j\in[k-1]}
\stackrel{\iid}{\sim} F$ and $\mathscr Y_n = \left\{ Y_j
\right\}_{j\in[n]} \stackrel{\iid}{\sim} G$, two-sample tests query the
hypotheses
\begin{align}
    &\mathcal H_0: F = G, &\mathcal H_1: 
    F \neq G.
    \label{eq:two_sample_test}
\end{align}
These tests have been studied
extensively in both theoretical \autocite[e.g.,][and references
therein]{bonnini2014nonparametric, thas2010comparing} and applied
\autocite[][]{charmpi2015weighted,conradsen2003test,schepsmeier2019goodness}
contexts, and see widespread application in science and industry.

Two-sample testing is well understood in the following two extremes:
\begin{enumerate}[label=(\alph*)]
    \item If $F = F^*$ is fully specified and $G \in \{ F^*, G^* \}$ is known to take on only a single alternative distribution, then the likelihood-ratio test
    (which ignores the $X$ samples) is optimal for fixed size.  
    \item If $F$ and $G$ can be arbitrary, then typically no
    best test exists under most notions of optimality, and general
    non-parametric tests based on, e.g., empirical CDFs 
    \autocite[for example,][]{kolmogorov1933sulla} are popular choices.
\end{enumerate}
In practice, one often encounters combinations of these two scenarios,
where likelihood-ratio type statistics are appealing but difficult to
control. For instance, if $F$ and $G$ are
known to ``cluster'' around two specified distributions $F^*$ and $G^*$, respectively
(e.g., the user has priors $\nu_F$ and $\nu_G$ on the space of probability
measures, whose expectations are $F^*$ and $G^*$), the likelihood ratio of
$F^*$ and $G^*$ has attractive power properties (it maximizes true
positives averaged over $\nu_G$), but difficult to control size (it only
controls false positives at a fixed rate averaged over $\nu_F$). The need
for these types of semi-parametric hypothesis tests arises naturally when
the data generating mechanism is broadly understood, but specific details
remain opaque. This situation arises frequently in modern science; for example, 
a practitioner might understand the biological principles underlying their
dataset well, yet may not have fully quantified the impact of measurement
noise (see e.g.~the discussion in \autocite{gao2005nonparametric}).
Currently, it is common practice to entirely forsake likelihood-type
approaches in such cases, and resort to the general non-parametric tests 
as in (b), trading desirable power properties for rigorous false-positive control.

Rank-based two-sample tests have emerged as a suitable tool to reconcile these two
divergent goals \autocite[see][and references therein]{gibbons2014nonparametric,klotz1962nonparametric}, providing efficient yet fully
distribution-free tests. Concretely, with $F_k$ and $H_{n,k}$ denoting the empirical
distributions of $\mathscr X_k$ and $\mathscr X_k\cup \mathscr Y_n$, it follows from
\autocite{chernoff1958asymptotic}, that, under suitable assumptions,
statistics of the form 
\begin{equation}
    T_{n,k}^J = \int J(H_{n,k}(x)) \ \mathrm{d}F_k(x) = \sum_{j=1}^k
    J(H_{n,k}(X_j))
    \label{eq:T_nk}
\end{equation}
are distribution-free, asymptotically normal as $n,k\to\infty, k/n \to \alpha>0$, and efficient against local alternatives $G$ for
a suitable choice of weight function $J = J_G$. In the case of location alternatives
$G(x) = G_n(x) = F(x-\mu/\sqrt{n})$, the test statistics resulting from appropriately chosen $J_G$ are the
popular Mann-Whitney $U$ \autocite{mann1947test} if $F$ is the logistic distribution, and
the Gaussian score transformed Mann-Whitney \autocite{van1956computation} if $F$ is
Gaussian. Moreover, \autocite{hodges1956efficiency} and \autocite{chernoff1958asymptotic}
showed that, in addition to performing favorably under logistic and Gaussian
$F$, the asymptotic efficiencies of these tests relative to the $t$-test are never below $\approx 0.86$ and
$1$, respectively, under \emph{any} $F$. These encouraging results prompted similar investigations in the context
of scale-alternatives $G_n(x)=F\left((1+\sigma/\sqrt{n}) x\right)$, where corresponding
choices of $J_G$ give rise to the Mood test \autocite{mood1954asymptotic}, Siegel-Tukey test
\autocite{siegel1960nonparametric}, and Gaussian score test \autocite{klotz1962nonparametric}.

Rank-based tests are increasingly used in small-sample settings \autocite[see][and references
therein]{mollan2020precise}, where their favorable power properties have been
confirmed to persist empirically. However, due to the slow convergence of their
associated central limit theorems, controlling the size of these tests in non-asymptotic settings is often difficult 
\autocite{conover1981comparative}, and there is a need for alternative methods of characterizing the finite-sample null distributions of rank-based test statistics. 

One of the contributions of this paper is to achieve this for
a closely related, asymptotically equivalent family of statistics based on rank spacings, which we now describe.
Let $X^{(j)}$ be the $j^{\text{th}}$ order statistic of $\mathscr X_k$, with conventions
$X^{(0)} = -\infty$ and $X^{(k)}=+\infty$ (and $F_k, H_{n,k}$ adjusted accordingly).
\autocite{holst1980asymptotic} showed that statistics of the form
\begin{equation*}
    Q_{n,k} = \sum_{j=1}^{k} w(F_k(X^{(j-1)}))\left( H_{n,k}(X^{(j)}) -
    H_{n,k}(X^{(j-1)})\right)
\end{equation*}
are asymptotically equivalent to 
$T_{n,k}^J$ in \eqref{eq:T_nk} when $w=J$. The difference $H_{n,k}(X^{(j)}) - H_{n,k}(X^{(j-1)})$ is
called a \emph{rank spacing}. Collecting these into a vector gives the equivalent representation
\begin{equation}
\label{eq:Snk1}
    \tilde{Q}_{n,k} = \sum_{j=1}^{k} w\left(\frac{j-1}{k}\right) S_{n,k}(j) = 
    \| \bfS_{n,k} \|_{1,w},
\end{equation}
where $(n+k)\bfS_{n,k}\in\mathbb Z_{\geq 0}^{k}$ with components $S_{n,k}(j) =
H_{n,k}(X^{(j)}) - H_{n,k}(X^{(j-1)}) - \frac{1}{n+k}$. Assuming continuous $F,G$ for the
moment, the additional $(n+k)^{-1}$ term allows for the convenient interpretation of
$S_{n,k}(j)$ as
\begin{equation*}
    (n+k)S_{n,k}(j) = \#\left\{ m : X^{(j-1)} < Y_m < X^{(j)} \right\},
\end{equation*}
and evidently does not alter the power of $Q_{n,k}$. The statistics $T_{n,k}^J$ and $\|
\bfS_{n,k} \|_{1,w}$ generally are not equivalent for finite samples (though they are in certain cases, e.g., Mann-Whitney's~$U$), but we will show that their power
properties are comparable for most statistical purposes. In \secref{sec:p=1}, we characterize the distribution
of $\| \bfS_{n,k} \|_{1,w}$ for arbitrary $n,k$, thereby enabling control of the size of tests based on this family of statistics in a precise way.

The $\bfS_{n,k}$ representation above naturally suggests a broader family of test
statistics $\left\{ \| \bfS_{n,k} \|_{p,w}^p \right\}_{p\geq 1}$ obtained by replacing the $1$-norm in \eqref{eq:Snk1} by the $p$-norm in the obvious way.
Such statistics arise in non-i.i.d.~two-sample testing (see Example 3 in the Supplementary Material) and in various applied
contexts \autocite{palamara2018high, riley2007locational}. The case $p=2,
w\equiv 1$, known as \emph{Dixon's statistic}, has received particular attention for its optimality properties in the context
of circular data \autocite{dixon1940criterion,weiss1956certain,gatto2015two,sethuraman1970pitman}; it is 
also connected to Greenwood's statistic  \autocite{greenwood1946statistical} in the limit of $k\to\infty$ and
$n$ fixed. Understanding the distributional
properties of the latter has been the subject of extensive study \autocite{moran1947random, moran1951random, moran1953random, gardner1952greenwood, darling1953class,burrows1979selected, currie1981further, stephens1981further, schechtman2000concentration}, yet a satisfactory description of its right-tail behavior
(which typically is the one of interest in testing goodness-of-fit) for finite samples
has remained elusive. In \secref{sec:p>1}, we fill this gap by characterizing the moments of
$\| \bfS_{n,k} \|_{p,w}^p$, and use this information to compute its CDF near the right boundary of its support. Additionally, we devise an algorithm to reconstruct the distribution of $\| \bfS_{n,k}
\|_{p,w}^p$ to $\varepsilon$ accuracy in $\mathcal O(\frac{kn}{\varepsilon}\log
\frac{n}{\varepsilon})$ time, paving the way for computationally efficient hypothesis testing.

Given that non-parametric statistics can match the efficiency of
likelihood-ratios in simple two-sample tests, while being exact for finite
samples, it is desirable to extend such a framework to the setting of
composite alternatives, where the relevant comparison is to the generalized
likelihood ratio test (gLRT). In Section \ref{sec:|A|>1}, we show that, for scale families, choices of $w$ mirroring the
Mann-Whitney and Gaussian score transformed Mann-Whitney test  \autocite{siegel1960nonparametric,
klotz1962nonparametric} do not exhibit similarly favorable power properties
as in the location setting. 
This confirms a conjecture of
\autocite{klotz1962nonparametric}, and 
suggests combining distinct
weight choices in a manner analogous to the gLRT. Using the
moment-reconstruction algorithm described above, we develop such a
technique in both the finite-sample and asymptotic regimes, and demonstrate empirically that the resulting tests can be powerful compared even to the gLRT.

Proofs of all the results presented are given in the Appendix.

\section{The case $p=1$}
\label{sec:p=1}

This section develops tools to compute the exact distribution of 
$\|\bfS_{n,k}\|_{1,w}$ defined in \eqref{eq:Snk1}. Assuming that $w\in \mathcal C^{2}[0,1]$, and expanding it
appropriately demonstrates that
\begin{align*}
    \| (n+k) \bfS_{n,k} \|_{1,w} = R_{n,k}^w  + c_S + \frac{1}{2}\varepsilon_S \quad \text{and}\quad T_{n,k}^w = R_{n,k}^w  + c_T + \frac{1}{2}\varepsilon_T,
\end{align*}
where $R_{n,k}^w = \sum_{j=1}^k H_{n,k}\left(X^{(j)}\right)w'\left( \frac{j-1}{k} \right)$;  $c_S, c_T$ are constants depending only on $n,k$ and $w$; and
\begin{align*}
    &\varepsilon_S = \sum_{j=1}^k H_{n,k}\left( X^{(j)} \right)
    \int_{\frac{j-1}{k}}^{\frac{j}{k}} kw''(x)\left( \frac{j}{k} - x \right) \ \mathrm dx
    \\
    &\varepsilon_T = \sum_{j=1}^k \int_{\frac{j-1}{k}}^{H_{n,k}\left( X^{(j)} \right)}
    w''(x)\left( H_{n,k}\left( X^{(j)} \right) - x \right) \ \mathrm dx.
\end{align*}

For $G$ sufficiently close to $F$ (or $w''$ appropriately small), these
error terms $\varepsilon_S, \varepsilon_T$ are generally $O(1)$ compared to
the $O(k)$ order of $R_{n,k}^w$, explaining the asymptotic equivalence of
$\| \bfS_{n,k} \|_{1,w}$ and $T_{n,k}^w$. Moreover, their similarity suggests
that even in finite-sample regimes, $\| \bfS_{n,k} \|_{1,w}$ and $T_{n,k}^w$
should generally behave comparably as long as $w$ is regular enough. We
do not quantify this statement precisely, but demonstrate that it is borne
out empirically in simulation studies like the one given in Supplementary
Figure~\ref{fig:roc}, where it is shown that, for  fixed $n=10,k=5$, the ROC curves of $\| \bfS_{n,k} \|_{1,w}$ and $T_{n,k}^w$ for various choices of $F, G$, $w$ (including highly irregular
ones)  match each other closely, indicating
that the favorable power properties of $T_{n,k}^w$
\autocite{conover1981comparative} are expected to transfer to $\| \bfS_{n,k}
\|_{1,w}$ as well. Therefore, it is of interest to study the finite-sample
distribution of $\|\bfS_{n,k}\|_{1,w}$.

The main result of this section is the following characterization of the law of
$\|\bfS_{n,k}\|_{1,w}$.  In what follows, we define $w_j=w(\frac{j-1}{k})$ and write $\bfw\in\mathbb R^{k}$ to denote $(w_1,\ldots,w_{k})$.

\begin{theorem}\label{thm:explicit_discrete}
    Let $\bfw\in\mathbb R^{k}$  have pairwise distinct entries and $w_{\max} = \max_{j=1}^{k} |w_j|$. Then the Laplace transform of $\|\bfS_{n,k}\|_{1,\bfw}$ is given by
    \begin{multline}
    \hspace{-2mm}    \mathbb E e^{t\|\bfS_{n,k}\|_{1,\bfw}} = (k-1)(-1)^{k+1}\cdot e^{t n w_{\max}} \times 
        \sum_{j=1}^{k-1} a_j^{e^{t(w_j-w_{\max})}} \bigg[ b_{n,k}
        \left( 1 - e^{t(w_j - w_{\max})(n+k-1)}\right)
        \\+ \sum_{m=0}^{k-3} c_{n,k,m} \left( 1 - e^{t(w_j - w_{\max})}
        \right)^{k-2-m} \bigg],
        \label{eq:explicit_discrete}
    \end{multline}
    where for any $\bfr\in\mathbb R^{k}$, $a_j^{\bfr} = \prod_{m\neq j} (r_j-r_m)^{-1}$ and  
    \begin{align*}
        &b_{n,k} = \frac{(-1)^{k}}{n+k-1}\cdot\binom{n+k-2}{k-2}^{-1},
        &c_{n,k,m} = \frac{(-1)^m}{n+1}\cdot
        \frac{\binom{k-2}{m}}{\binom{n+m+1}{m}}.
    \end{align*}
\end{theorem}

\begin{remark}
For hypothesis testing, \eqref{eq:explicit_discrete} needs to
be inverted in order to recover the requisite null distribution.
This can be done quickly and in a numerically stable manner, as we demonstrate in Supplementary Figure~\ref{fig:cdf}, where CDFs obtained from Monte Carlo
iterates are contrasted with those computed from \eqref{eq:explicit_discrete}.
\end{remark}

\begin{remark}\label{rmk:discrete}
    The assumption that the components of $\bfw$ are distinct is merely to simplify equation \eqref{eq:explicit_discrete}, and can be dropped. In case there are ties, the result
    is
    obtained by evaluating \eqref{eq:explicit_discrete} along a sequence $\bfw^{(n)}$ of
    weights whose entries
    are mutually distinct and converge to $\bfw$. Explicit expressions (involving
    suitable partial
    derivatives of $a_j^{\bfw}$ in each component of $\bfw$) can be found in the Supplementary
    Material. Numerical
    inversions are performed without difficulty as before.
\end{remark}

Theorem \ref{thm:explicit_discrete} enables hypothesis testing in regimes
of $n$ and $k$ both remaining small, thereby complementing results of
\autocite{holst1980asymptotic} where the asymptotic behavior of $\| \bfS_{n,k}
\|_{1,\bfw}$ for $n,k\to\infty$, $k/n\to\alpha$ is considered. This leaves
open the case of one parameter, say (without loss of generality) $n$,
diverging towards $\infty$, with the other, $k$, kept fixed. With
experimental methods producing ever more refined, yet possibly sparse,
data, this situation is encountered increasingly often 
\autocite[see, e.g.,][for perspectives from biology and engineering]{he2009learning,yang2006impact}. The following result characterizes $\|
\bfS_{n,k} \|_{1,\bfw}$ in this regime.

\begin{theorem}\label{thm:explicit_continuous}
    As $n\to\infty$ with $k$ remaining fixed, $\mathbb P\left( \| \bfS_{n,k}
    \|_{1,\bfw} \leq x\right) = \mathbb P\left( \| \bfS_k \|_{1,\bfw} \leq x \right) +
    \varepsilon(x)$, where $\bfS_k \sim
    \operatorname{Dirichlet}(\mathbf{1}_{k})$
    (with $\mathbf{1}_{k+1}\in\mathbb Z^{k}$ being the all-ones vector)
    is uniformly distributed on the $(k-1)$-dimensional simplex, and $\|
    \varepsilon\|_{\infty}\in O(n^{-1})$. Moreover, with $a_j^{\bfw}$ as in Theorem~\ref{thm:explicit_discrete},
    \begin{equation}
        \mathbb P\left( \| \bfS_k \|_{1,\bfw} \leq x \right) = (-1)^{k+1} \sum_{j=1}^{k-1} a_{j}^{\bfw}\left[ (x-w_j)_+
        \right]^{k-1},
    \label{eq:explicit_continuous}
    \end{equation}
    as long as the components of $w$ are pairwise distinct.
\end{theorem}

\begin{remark}
    As with Theorem \ref{thm:explicit_discrete}, the distinctness assumption on $\bfw$ can
    be relaxed by taking suitable limits.
\end{remark}

\section{The case $p>1$}
\label{sec:p>1}

The explicit form of Theorem \ref{thm:explicit_discrete} relies on the observation that
$\bfS_{n,k} \sim \operatorname{Multinomial}(n, \bfS_k)$ allows factorization of the Laplace
transform of $\| \bfS_{n,k} \|_{1,\bfw}$: $\mathbb E e^{t\| \bfS_{n,k} \|_{1,\bfw}} = \mathbb E \|
\bfS_k \|_{1,e^{t\bfw}}^n$. When $p > 1$, interaction terms in $\| \bfS_{n,k} \|_{p,\bfw}^p$ prevent such a factorization. Nevertheless, the individual moments can still be accessed.

\begin{theorem}
    Let $G_p(x,y) = \sum_{m=0}^{\infty} \mathrm{Li}_{-pm}(x)y^m/m!$, where
    $\mathrm{Li}_s(x) = \sum_{j=1}^{\infty} j^{-s} x^j$ is the polylogarithm function.
    Denoting by $[x^ny^m] P(x,y)$ the $(n,m)^{\text{th}}$ coefficient of a formal power
    series $P$ in $x$ and $y$, we have
    \begin{equation}
        \mathbb E \left(\|\bfS_{n,k}\|_{p,\bfw}^{p}\right)^m =
        \dfrac{m!}{\binom{n+k-1}{k-1}} [x^ny^m] \prod_{i=1}^k G_p\left( x, w_i y
        \right).
        \label{eq:explicit_discrete_moments}
    \end{equation}
    In particular, the first $m$ moments of $\|\bfS_{i,j}\|_{p,\bfw}^p$ for $(i,j) \in
    \{0, \dots, n\} \times \{1, \dots, k\}$ can be computed in
    $O\left(nm\cdot  (\log nm) \cdot k\right)$ time.
    \label{thm:explicit_discrete_moments}
\end{theorem}

As with the $p=1$ case, there are three regimes of interest:
\begin{enumerate}
    \item $n,k\to\infty$ while $k/n \to \alpha$; \label{item:r1}
    \item  $n,k$ both small; and\label{item:r2}
    \item $n\to\infty$ with $k$ fixed.\label{item:r3}
\end{enumerate}
Regime \ref{item:r1} is covered by the same central limit theorems in
\autocite{holst1980asymptotic} that resolved the corresponding question when $p=1$.
Theorem \ref{thm:explicit_discrete_moments} will turn out to be
useful primarily in regime \ref{item:r2}, while the following analogue of Theorem~\ref{thm:explicit_continuous} covers regime \ref{item:r3}.

\begin{theorem}\label{thm:explicit_continuous_moments}
    Let $Q_p(x) = \sum_{m=0}^{\infty} (pm)!x^m/m!$.  Then for $n\to\infty$ with $k$ kept
    fixed, $\mathbb P(\| \bfS_{n,k} \|_{p,\bfw}^p \leq x) = \mathbb P\left( \| \bfS_k
    \|_{p,\bfw}^p \leq x \right) + \varepsilon(x)$, where $\| \varepsilon \|_{\infty} \in
    O(n^{-1})$, and
    \begin{equation}
        \mathbb E\left( \| \bfS_k \|_{p,\bfw}^p \right)^m =
        \dfrac{(k-1)!m!}{(pm+k-1)!} [x^m] \prod_{j=1}^k Q_p(w_jx).
        \label{eq:explicit_continuous_moments}
    \end{equation}
    In particular, the first $m$ moments of $\| \bfS_j \|_{p,\bfw}^p$ for $j\in\{1,
    \dots, k\}$ can be computed in\break $O(m\cdot (\log m)\cdot k)$ time.
\end{theorem}

\begin{remark}
The generating function $Q_p(x)$ can be expressed as the generalized hypergeometric
series $Q_p(x) = {}_pF_0\big[ 1, \frac{1}{p},
\frac{2}{p}, \dots, \frac{p-1}{p} \big]( p^2x).$
In particular, for $p=2$ (i.e., including the Greenwood statistic 
$\|\bfS_k\|_{2,\bfmath{1}_k}^2$), we have $Q_2(x) = {}_2F_0\big[ 1,\frac{1}{2} \big](4x)
=\frac{1}{\sqrt{x}}D\big( \frac{1}{2\sqrt{x}} \big)$, where Dawson's
integral
\begin{equation*}
    D(x) = e^{-x^2}\int_0^x e^{t^2} ~\mathrm{d}t
    \label{eq:dawson}
\end{equation*}
is interpreted through its asymptotic expansion 
\autocite[cf.][formula 7.1.23]{abramowitz1965handbook}.
\end{remark}

In order to perform hypothesis testing, Theorems
\ref{thm:explicit_discrete_moments} and \ref{thm:explicit_continuous_moments} require
numerical ``inversion'' similar to the need for inverse Laplace transforms to render
\eqref{eq:explicit_discrete} and \eqref{eq:explicit_continuous} practical. More
concretely, they require efficient reconstruction of a (compactly supported) distribution
$F$ given its truncated moment sequence $\int x^j \ dF(x), j\in[m]$.
Although this is a well-studied problem 
\autocite[see][for a comprehensive introduction]{akhiezer2020classical}, equations 
\eqref{eq:explicit_discrete_moments} and
\eqref{eq:explicit_continuous_moments} have some unique properties that are not often encountered:

\begin{enumerate}[label=(\alph*)]
    \item An arbitrary number of moments can be computed efficiently. This
        is markedly distinct from situations in which moments are estimated
        from, e.g., experimental observations and limited in number.
        Applying tools developed in this latter context \autocite[ranging from extremal inequalities to maximum entropy based approaches to concretely applications-driven methods, see, e.g.,][]{schmudgen2020ten, john2007techniques} typically
        under-utilizes all the information available, or becomes computationally
        infeasible.
    \item Each moment can be computed exactly, and therefore usual concerns around the
        well-conditioning of the moment problem \autocite[see, for example,][]{talenti1987recovering} do not apply.
\end{enumerate}

By exploiting these two properties, we can carry out the necessary reconstruction efficiently and with great accuracy, simply by considering expectations of Bernstein polynomials.

\begin{proposition}
    Let $X\in[0,1]$ be a random variable either (a) continuous with density $f\in C^1\left([0,1]\right)$,
    or (b) discrete with support $\mathrm{supp}_X = \{x_0, \dots, x_N\}$, and $F$ be its CDF. Moreover,
    denote by $B_{n,x}$ the degree-$n$ Bernstein polynomial approximating
    $\mathbbm{1}_{[0,x]}$. Then, for any resolution $\varepsilon_n \to 0, 
    \varepsilon_n > n^{-1/2}$, there exists $n_0(f,\varepsilon) \in 
    \mathbb N$, so that for all $n\geq n_0$,
    \begin{align}
        \sup_{x\in[0,1]} \big| \mathbb EB_{n,x}(X) - F(x) \big| &\leq
        \frac{\|f\|_{\infty} + 2\|f'\|_{\infty}+2}{n+1} \tag{a},
        \label{eq:bernstein_convergence_i} \\
        \sup_{x\in [0,1]\setminus \mathrm{supp}_X^{\varepsilon_n}} \big
        | \mathbb E B_{n,x}(X) - F(x) \big| &\leq e^{-2n\varepsilon_n^2}, 
        \tag{b} \label{eq:bernstein_convergence_ii}
    \end{align}
    where $\mathrm{supp}_X^{\varepsilon} = \left\{ x\in[0,1] ~:~ d(x, \mathrm{supp}_X)
    < \varepsilon \right\}$ is the $\varepsilon$-fattening of $\mathrm{supp}_X$.
    \label{prop:hausdorff}
\end{proposition}

Several features of the proposition are worth highlighting:
\begin{enumerate}
    \item By virtue of $B_{n,x}$ being a degree-$n$ polynomial, $\mathbb EB_{n,x}(X)$ is
        just a linear combination of the first $n$ moments $\mu_1, \ldots, \mu_n$ of $X$;
        more explicitly,
        \[
            \mathbb EB_{n,x}(X) = \sum_{m=0}^{\lfloor nx \rfloor}\binom{n}{m}(-1)^{n-m}
            \left(\delta^{n-m}\mu\right)_m,
        \]
        where $\mu = \left( \mu_j \right)_{j\in \mathbb N}$ denotes the moment sequence
        of $X$, and $\delta$ is the difference operator.
    \item For a discrete $X$ of $n$ atoms, $n$ moments are sufficient to determine the
        distribution of $X$ via solving an $n\times n$ Vandermonde system, which
        can be performed in $O(n^2)$ time \autocite{bjorck1970solution}. However,
        $\|\bfS_{n,k}\|_{p,\bfw}^p$ generically has $O(n^{k-1})$ atoms, whose precise
        locations within $\{x_{\min}, \dots, \|\bfw\|_{\infty} n^p\}$ are typically
        unknown, therefore requiring $O( \min\{\|\bfw\|_{\infty}^2n^{2p}, 
        n^{2(k-1)} \} )$ operations, which is prohibitively large even for small values
        of $p$ or $k$.
    \item $B_{n,x}$ may be replaced with any other polynomial approximation scheme in
        order to impose desired properties on the reconstructed density. For instance, if the user wishes to perform a one-sided test, then resorting to one-sided polynomial 
        approximations \autocite[the optimal of which is worked out in][]{bustamante2012best} is more suitable.
\end{enumerate}

Beyond its practical impact in performing two-sample tests when $n$ is large and $k$
modest, the quantity $\| \bfS_k \|_{p,\bfw}^p$ appearing in Theorems
\ref{thm:explicit_continuous} and \ref{thm:explicit_continuous_moments} is of 
independent interest in the context of one-sample testing, where it constitutes the
appropriate equivalent of $\| \bfS_{n,k} \|_{p,\bfw}^p$. The case of Greenwood's statistic
(corresponding to $p=2$ and $\bfw = \mathbf{1}_{k}$) has received particular attention, 
with extensive studies clarifying left-tail behaviour, asymptotic normality as 
$k\to\infty$, and large deviation functions. Theorem
\ref{thm:explicit_continuous_moments} can be used to supplement these results with a
characterization of the right tail.

\begin{proposition} \label{prop:right_tail}
    Without loss of generality, assume $\bfw\in\mathbb R_+^{k+1}$ and $\|\bfw\|_{\infty} =
    1$, and denote by
    \[
    W_{\bfw} = \left| \{ 1 \leq j \leq k+1 ~ : ~ w_j = 1 \} \right|
    \]
    the
    number of weight components assuming value $1$. Then the density $f_k^{p,\bfw}$ of 
    $\| \bfS_k \|_{p,\bfw}^p$ is analytic on $[x_0,1]$, where
    \[
        x_0 = 
        \begin{cases}
            \frac{1}{2^{p-1}}, &\text{if\,  $W_{\bfw} = k+1$}, \\
            \max_{j: w_j < 1} w_j, &\text{otherwise},
        \end{cases}
    \]
    and its degree-$r$ Taylor polynomial around $1$ can be computed in $O\left(
    \frac{r}{p}\log\frac{r}{p} \log k + [r\log r]^2\right)$ time. For $r=k-2$ it reads
    \begin{equation*}
        f_k^{p,\bfw}(x) = \dfrac{(k-1)W_{\bfw}}{2^{k-1}}\left( 1-x \right)^{k-2} +
        O\left( (1-x)^{k-1} \right).
        \label{eq:general_taylor}
    \end{equation*}
    In particular, Greenwood's statistic satisfies
    \begin{equation*}
        f_k^{2,\mathbf{1}_{k}}(x) = \dfrac{\binom{k}{2}}{2^{k-2}} 
        \left( 1-x \right)^{k-2} + O\left( (1-x)^{k-1} \right).
        \label{eq:taylor_expansion}
    \end{equation*}
\end{proposition}

The right tail is typically the one of interest in one- and two-sample tests, and
so as long as long as the desired significance threshold $\alpha$ is less than 
$\mathbb P\left( \| \bfS_k \|_{p,\bfw}^p \geq x_0 \right)$, Proposition \ref{prop:right_tail} 
allows for calculating $\varepsilon$-accurate $p$-values in $O(\log^2\varepsilon)$ time.
This compares favorably with the $O\left(\varepsilon^{-1}\right)$ rate of Theorem~\ref{thm:explicit_continuous_moments},
and can provide a substantial speed-up for
large data sets.

\section{Hypothesis testing when $|\mathcal A| > 1$}
\label{sec:|A|>1}

Assume without loss of generality that $X\sim \operatorname{Uniform}([0,1])$ and 
$Y$ has density $g(x) = 1 + h(x)/\sqrt{n}$. In the case of singleton hypotheses
$F = F^*$ and $G \in \{ F^*, G^*\}$, 
$\| \bfS_{n,k} \|_{1,\bfh}$ can be regarded as a
non-parametric version of the likelihood ratio test for alternatives $G$ that are near
$F$. This follows from the asymptotic equivalence between 
tests based on $\| \bfS_{n,k} \|_{1,\bfw}$ and likelihood-ratio tests \autocite{holst1972asymptotic}, and 
can also be seen by observing that
\begin{multline*}
    \frac{\sqrt{n}}{n+k}\log\prod_{j=1}^n g(Y_j) = \frac{\sqrt{n}}{n+k}
    \sum_{j=1}^n \log \left[ 1 + \frac{h(Y_j)}{\sqrt{n}}\right]
    \approx \frac{1}{n+k}\sum_{j=1}^n h(Y_j) \approx \sum_{j=1}^{k}
    h\left(\frac{j-1}{k}\right) S_{n,k}(j) = \|\bfS_{n,k}
    \|_{1,\bfh},
\end{multline*}
where $\bfh \in\mathbb R^{k+1}$ has $j^{\text{th}}$ component
$h\left((j-1)/k\right)$.

By analogous reasoning, if the alternative hypothesis $G \in \big\{ 1+h^{\theta}(x)/\sqrt{n} \big\}_{\theta \in \Theta}$ is composite 
(and parameterized by $\theta$ over some index set $\Theta$), then given observations
$X_1, \ldots, X_{k-1}$ and $Y_1, \ldots, Y_n$, one may expect tests based on $\sup_{\theta\in
\Theta} \| \bfS_{n,k} \|_{1,\bfh^{\theta}}$ to provide non-parametric equivalents of
generalized likelihood-ratio tests. When $|\Theta| = m < \infty$,  multivariate
extensions of the previous results follow in a straightforward manner.

\begin{proposition}\label{prop:multivariate_extension}
    For $m$ weights $\bfw^1, \ldots, \bfw^m\in\mathbb R^{k}$, each with pairwise distinct entries,
    the Laplace transform of the tuple $S_{(m)} = \left( \|\bfS_{n,k}\|_{1,\bfw^1}, \ldots, \|
    \bfS_{n,k} \|_{1,\bfw^m} \right)$ is given by
    \begin{multline*}
        \mathbb E e^{\langle t, S_{(m)} \rangle} = (k-1)(-1)^{k+1}\cdot e^{n w^{\max}}
         \sum_{j=1}^{k-1} a_j^{e^{\omega_j} } \bigg[ b_{n,k} \left( 1 - e^{\omega_j
        (n+k-2)} \right) + \sum_{m=0}^{k-3} c_{n,k,m} \left( 1 - e^{\omega_j}
        \right)^{k-1-m} \bigg],
        \label{eq:multivariate_extension}
    \end{multline*}
    where $t=(t_1,\ldots,t_m)$, $\omega_j = \sum_{r=1}^m t_r w_j^r - w^{\max}$ with $w^{\max} = \max_j
    \sum_{r=1}^m t_r w^r_j$, and $a_j^w, b_{n,k}$ and $c_{n,k,m}$ are defined as in Theorem~\ref{thm:explicit_discrete}.
    
    Moreover, the joint moments of $S_{(m)}$ can be computed in $O\big(n\prod_{j=1}^r
    m_j$ $\times  (\log n\prod_{j=1}^r m_j) \times k\big)$ time as
    \begin{equation*}
        \mathbb E\prod_{j=1}^r \left(\|\bfS_{n,k}\|_{p_j,\bfw^j}^{p_j}\right)^{m_j} =
        \dfrac{\prod_{j=1}^r m_j!}{\binom{n+k-1}{k-1}} [x^ny_1^{m_1}\cdots y_r^{m_r}]
        \prod_{i=1}^k
        G_r\left( x, w^1_i y_1, \dots, w^r_i y_r \right),
        \label{eq:main_mulitvariate}
    \end{equation*}
    where $G_r(x,y_1, \dots, y_r) = \sum_{m_1, \dots, m_r=0}^{\infty}
    \mathrm{Li}_{-\sum_{j=1}^r p_jm_j}(x)\prod_{j=1}^r y_j^m/m_j!$. These joint moments
    can be used to approximate $\mathbb P\left( \| S_{(m)} \|_{\infty} \leq x \right)$
    up to $\varepsilon$ accuracy in $O\left( \varepsilon^{-1} \right)$ time.
\end{proposition}

Part of the motivation for formulating Proposition \ref{prop:multivariate_extension} is 
to improve the performance of non-parametric testing procedures in the context of scale
alternatives. As noted earlier, for location families, the weight functions $w_1^{\mu}(x) = x$ and $w_2^{\mu}(x)
= \Phi^{-1}(x)$, where $\Phi$ denotes the standard Gaussian CDF, are known to compare impressively
against the parametric $t$-test when alternatives $G_n$ are shifts of $F$, with Pitman efficiencies never dropping below $\approx 0.86$ and $1$, respectively \autocite{hodges1956efficiency}.
However, for scale families, the corresponding choices 
$w_1^{\sigma}(x) = (x-1/2)^2$ and $w_2^{\sigma}(x) = \Phi^{-1}(x)^2$ \autocite{ansari1960rank} compare less favorably against the relevant $F$-test:
\cite{sukhatme1957certain} demonstrated that $\inf_{F} e_{w_1^{\sigma}}^{\sigma}(F) = 0$ (where
$e_{w_1^{\sigma}}^{\sigma}(F)$ denotes the Pitman efficiency of $\| \bfS_{n,k}
\|_{1,\bfw_1^{\sigma}}$ against the $F$-test under scale shifts) under $X\sim F$, while
\autocite{klotz1962nonparametric} showed that $\inf_F e_{w_2^{\sigma}}^{\sigma}(F) \leq
0.47$ and conjectured that, in fact, efficiencies arbitrarily close to zero can be
realized. The following example confirms a considerably stronger version of Klotz'
conjecture.

\begin{proposition}\label{prop:klotz}
    For $a,b>0$, define random variables $X_{a,b}$ through their densities
    \[
        f_{a,b}(x) = \frac{1}{Z} \times
            \begin{cases}
                 \frac{1}{a} \frac{x+b}{b-a},    & \text{if }-b \leq x < -a, \\
                 -\frac{1}{x},                   & \text{if }-a \leq x < -1, \\
                 1,                              & \text{if }-1 \leq x < 1,  \\
                 \frac{1}{x},                    & \text{if }1  \leq x < a,  \\
                 \frac{1}{a} \frac{x-b}{a-b},    & \text{if }a  \leq x < b,
            \end{cases}
    \]
    where $Z = 2\log a + 1 + \kappa$, with $\kappa = b/a$. Then, as $a\to\infty$ while
    keeping $\kappa\in o(\log a)$, $e_w^{\sigma}(F_{a,b})\to 0$ for any $w\in
    \mathcal C^1([0,1])$ whose derivative is bounded by $C [ x^{-1}|\log x|^p +\break
    (1-x)^{-1}|\log(1-x)|^p]$ for some constants $C>0$ and $p>0$.
\end{proposition}

If it is unknown whether the data-generating mechanism for $X_1, \ldots, X_k$ lies within 
distributions against which weights like $w_2^{\sigma}$ are powerful  (which includes most 
named distributions, see \autocite{klotz1962nonparametric} and Example 1 in the Supplementary
Material), or is closer to the one
described in Proposition \ref{prop:klotz}, then combining $w_2^{\sigma}$ with a complementary
weight function in the manner outlined by Proposition \ref{prop:multivariate_extension}
may boost performance.

\section{Application to non-parametric hypothesis tests}
\label{sec:application}
\begin{figure}[t]
    \centering
    \includegraphics[width=\textwidth]{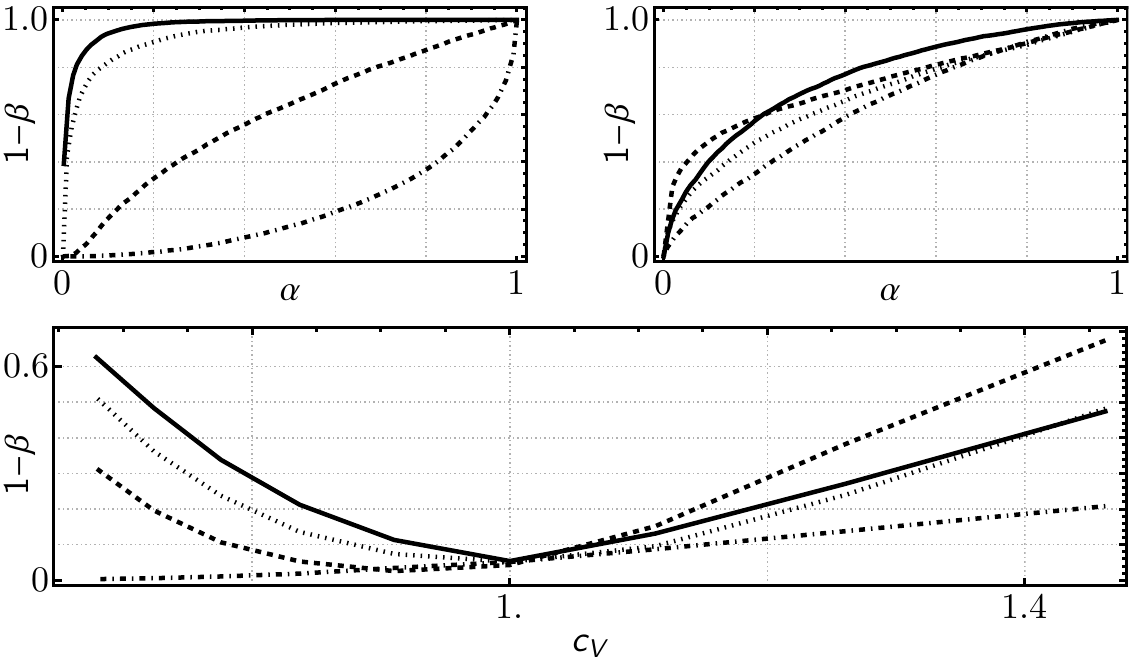}
    \caption{One-sample test comparison between Greenwood's
        $\|\bfS_k\|_{2,\mathbf{1}_{k}}^2$ (solid line), Cox and Oakes' $CO_k$ (dashed),
        Gail and Gastwirth' $G_k$ (dotted), and the Cram\'er-von Mises test
        (dot-dashed) for $k=20$. Upper panels display ROC curves of type-I error ($\alpha$)
        against power ($1-\beta$) in the case of under- and over-dispersed
        alternatives (left and right), respectively. Bottom panel illustrates
        power against varying coefficient of variation $c_V$.
    }
    \label{fig:greenwood}
\end{figure}
We begin by carrying out the original test of uniformity proposed by
\cite{greenwood1946statistical} for moderately sized $k=20$ (which does not yet induce
CLT-type behavior) and comparing it to three other common tests. This analysis extends
previous power studies that either omitted Greenwood's statistic for lack of exact
$z$-scores, or accepted approximation errors in their results \autocite[e.g.,][]{d1986goodness,henze2005recent}. Despite the small scale of our comparison, the
results are promising, and we hope they will encourage inclusion of Greenwood's statistic into
future benchmarking efforts.

\autocite{greenwood1946statistical} was interested in testing under- or
over-dispersion of spacings relative to a homogeneous Poisson Point Process;
that is, he considered the hypotheses
\begin{align*}
    &\mathcal H_0: \left\{ T_j \right\}_{j\in[k]} \stackrel{\iid}{\sim} 
    \text{Exp}(\lambda),
    &\mathcal H_1: \left\{ T_j \right\}_{j\in[k]} \stackrel{\iid}{\sim} X,
    \text{ where }c_{V}^2 = \frac{\mathrm{Var}X}{\left(\mathbb EX\right)^2} \neq
    1,
\end{align*}
which is equivalent to testing whether
$(T_1, T_2, \dots, T_k)/\sum_{j=1}^{k} T_j$ is
distributed like $\bfS_k$ under the null. Greenwood proposed to use
$\|\bfS_k\|_{2,\mathbf{1}_{k}}^2$, but
was not able to quantify its power, numerically or otherwise. Theorem
\ref{thm:explicit_continuous_moments} and Proposition \ref{prop:hausdorff} allow to
compute the law of Greenwood's statistics quickly (computing $p$-values of $10,000$
simulation runs each with $k\leq 30$ takes $\approx 5$ seconds on an ordinary
laptop), facilitating power comparisons. 

The test statistics we compared against  were $CO_k$ from \autocite{cox2018analysis}, $\bfS_k$ from
\autocite{gail1978scale}, and the Cram\'er-von Mises test.
\autocite{henze2005recent} identified these as high-performing goodness-of-fit tests through
extensive simulation studies. We compared them
to Greenwood's statistic on the same under- and over-dispersed alternatives
used in \citeauthor{henze2005recent}: the uniform distribution on $[0,1]$, and the Weibull
distribution of scale $1$ and shape $0.8$. The results, together with a sensitivity
analysis of power against varying dispersion (using the family of Weibull distributions
of scale $1$ and shapes $0.8,0.9,\dots,1.5$), are displayed in Figure
\ref{fig:greenwood}. They reveal competitive performance of
$\|\bfS_k\|_{2,\mathbf{1}_{k}}^2$, especially in the under-dispersed regime (upper-left).

\begin{figure}[h]
    \centering
    \includegraphics[width=\textwidth]{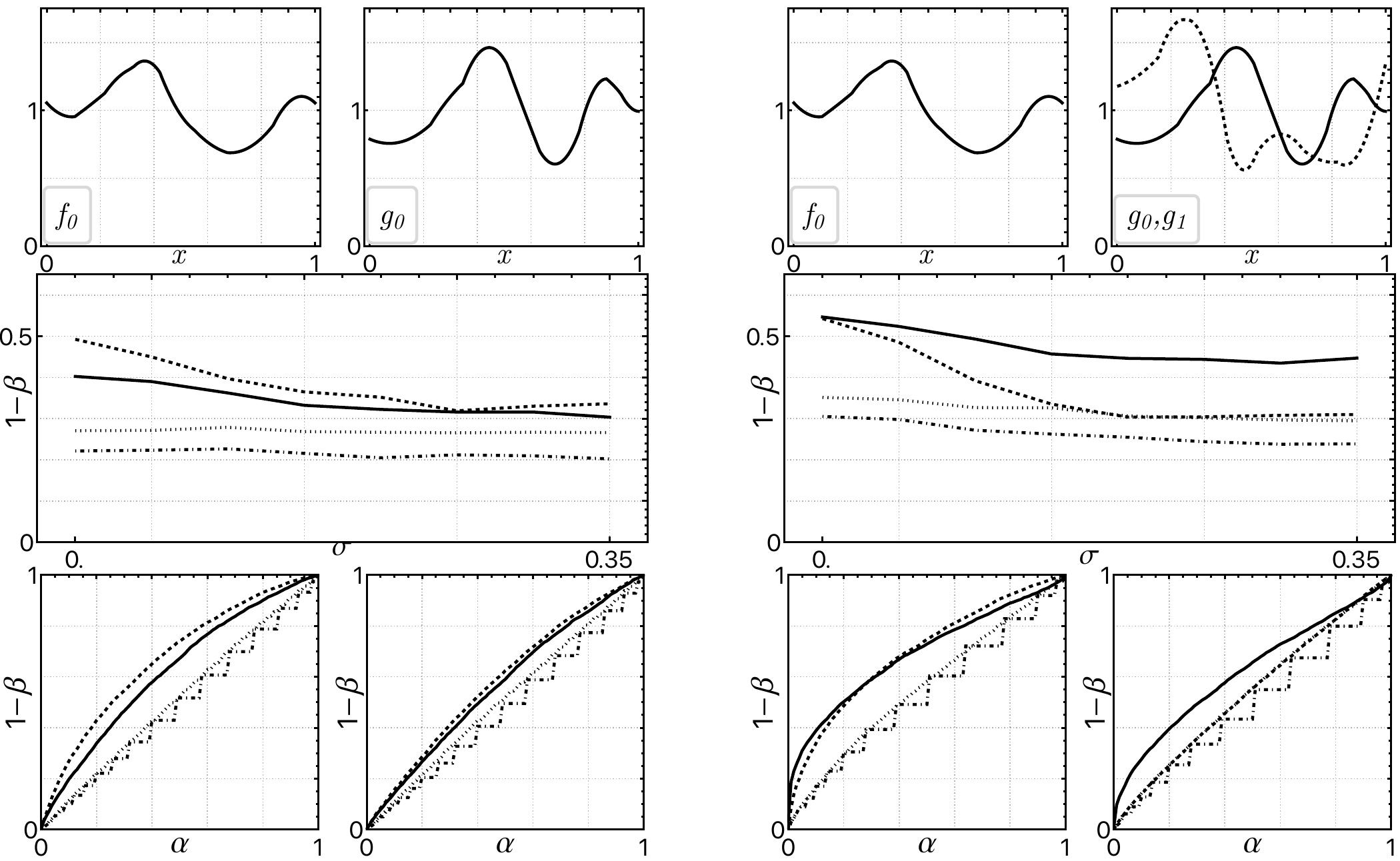}
    \caption{Two-sample test comparisons (center and bottom row) of $\| \bfS_{n,k}
        \|_{1,\bfh}$ (solid line), the likelihood-ratio statistic (dashed),
        Mann-Whitney's U statistic (dotted), and the Kolmogorov-Smirnov statistic
        (dot-dashed) on $\mathcal N$ and $\mathcal A$ centered around given distributions
        (top row). $(1-\beta)$ and $\alpha$ denote power and test size as in
        Figure~\ref{fig:greenwood}, and $\sigma$ the measurement noise. 
        Plots in the middle row correspond to $\alpha = 0.05$.
        In each scenario, ROC curves in the
        bottom row correspond to $\sigma=0$ and $0.3$, respectively. All
        simulations were run on $n=50, k=25$ and $10,000$ Monte-Carlo iterations.
}
    \label{fig:likelihood_ratio}
\end{figure}

To empirically probe the relevance of Theorem \ref{thm:explicit_discrete} and the 
multiple testing strategy presented in the previous section, we compared power properties
of $\| \bfS_{n,k} \|_{1,\bfh}$ and $\max\{ \| \bfS_{n,k} \|_{1,\bfh_0}, \| \bfS_{n,k} \|_{1,\bfh_1} \}$
against the gLRT, as well as
against two omnibus tests \autocite{kolmogorov1933sulla,mann1947test}
which are widely used in practice.
Results are
shown in Figure~\ref{fig:likelihood_ratio}, with simple and composite cases divided into
left and right column, respectively; and null ($f_0$) and alternative ($g_0=1+h_0/\sqrt{n},g_1 = 
1+h_1/\sqrt{n}$)
distributions were chosen to reflect fairly generic multi-modal two-sample setups (top
row). In order to simulate ``measurement noise'' or misspecification of $F$ and
$G$ around $f_0, g_0$ and $g_1$, $\mathscr X_k$ and $\mathscr Y_n$ samples were perturbed by Gaussian
noise of varying standard deviation $\sigma\in[0,0.35]$, and power for the (generalized)
likelihood-ratio statistic, its $\| \bfS_{n,k} \|_{1,\bfh}$ counterparts, and the two
omnibus tests was computed at size $\alpha = 0.05$ (center row of Figure
\ref{fig:likelihood_ratio}). As expected, the likelihood-ratio dominates in the noiseless
regime. However, $\| \bfS_{n,k} \|_{1,\bfh}$ and its generalized extension
perform competitively, closing the gap to likelihood ratio tests (or in the case of
composite alternatives, reversing it) as noise is introduced. Importantly,
calibrating likelihood-ratio tests in these contexts requires exact knowledge of the
perturbation (which in general is not accessible to the practitioner), while tests 
based on $\| \bfS_{n,k} \|_{1,\bfh}$ do not. A more thorough description of the various
compared tests for all sizes and fixed $\sigma = 0$ and $\sigma = 0.3$ is provided in the
ROC curves of Figure \ref{fig:likelihood_ratio} (bottom row), which confirm that the
favorable performance of $\| \bfS_{n,k} \|_{1,\bfh}$ persists across the range of
$\alpha$ most relevant in practice.

Although we find the theoretical and simulation evidence presented here convincing, this alone is not enough to ensure that our results will be utilized elsewhere. 
To aid practitioners in applying our methods,
we provide
code implementing most of the functionality outlined in this manuscript at
\url{https://github.com/songlab-cal/mochis} (currently as a Mathematica notebook, but
python and R packages are forthcoming). Its interface allows users to specify $f_0$ and
any number of $g_i \in \mathcal A$ on a bounded interval or $\mathbb R$ through either a
simple drag-and-drop mechanism or explicitly in closed or numerical form. From there,
the relevant distribution of $\|
\bfS_{n,k} \|^p_{p,\bfw}$ (in the case of $p=1$) or moments ($p>1$) are directly
computed, and $p$-values corresponding to a given set of samples $X_1, \ldots, X_k$, $Y_1,\ldots, Y_n$ calculated. Optional arguments allow customization of any part of the
procedure. Even though the current implementation focuses on the one- and two-sample
situations described above, several generalizations are straightforward to include:

\begin{enumerate} 
    \item When extending results from continuous variables to discrete ones, ties can be resolved
        uniformly at random when constructing $\bfS_{n,k}$ from $X_1, \ldots, X_k$ and $Y_1, \ldots, Y_n$.
    \item The i.i.d. assumption on $X_1, \ldots, X_k$ and $Y_1, \ldots, Y_n$ can be relaxed to any other
        setting where null distributions effectively reduce to uniform samples from the discrete or
        continuous simplex; e.g., the same reasoning applies to paired two-sample tests. 
    \item Several representative weight choices corresponding to commonly encountered alternatives
        (e.g., $\bfw = (k, k-1, \ldots, 0)$ associated with Mann-Whitney's $U$ statistic in the case of
        $Y$ stochastically dominating or being stochastically dominated by $X$) are included in the code
        base as pre-computed tables for the case of $p = 1$ due to their relevance in two-sample
        testing. An interface allows users to specify similar generic weight choices (not
        necessarily arising from any fixed $f_0$ and $g_0$) for both $p = 1$ and $p > 1$ (which
        can become relevant for non-i.i.d. data).
    \item The hypothesis testing results derived here only relied on the moments of $\| \bfS_k
        \|_{p,\bfw}^p$ and $\| \bfS_{n,k} \|_{p,\bfw}^p$ to reconstruct $\mathbb E 
        \left[\mathbbm{1}_{[0,t]}\left( \| \bfS_k \|_{p,\bfw}^p \right)\right]$ and $\mathbb E
        \left[
        \mathbbm{1}_{[0,t]}\left( \| \bfS_{n,k} \|_{p,\bfw}^p \right) \right]$ or their equivalents in the
        context of composite alternatives, where CDFs of maxima of $\| \bfS_k \|_{p,\bfw}^p$ and $\| \bfS_{n,k}
        \|_{p,\bfw}^p$ are of interest. Of course, $\mathbb E f\left( \| \bfS_k \|_{p,\bfw}^p
        \right)$ and $\mathbb
        E f\left( \| \bfS_{n,k} \|_{p,w}^p \right)$ can be approached in a similar fashion for any 
        $f\in L^2\left( [0,1] \right)$.
    \item The moment-reconstruction method described through Proposition \ref{prop:hausdorff} is
        applicable to any bounded random variable whose moments are known exactly, and can be used
        accordingly in the code implementation.
    \item Extension of the non-parametric generalized-likelihood-type test as formulated
        above to the asymptotic regime requires knowledge of the distribution of the
        maximum of an arbitrary number of correlated Gaussian variables, which in general
        is intractable. Switching to a simpler summary like the sum, however, is feasible
        and may offer similar power depending on the precise correlation structure.
        Analyzing the details of this situation is left for future work.
\end{enumerate}

\section*{Acknowledgments}  
We thank Ben Wormleighton for acquainting the authors with Ehrhart's work, and Jonathan Fischer for
helpful comments on software implementation. This research is supported in part by an NIH grant
R35-GM134922.  

\newpage

\printbibliography
\clearpage
\appendix

\setcounter{table}{0}
\renewcommand{\tablename}{{Table}}  
\renewcommand{\thetable}{{S\arabic{table}}}%
\setcounter{figure}{0}
\renewcommand{\figurename}{{Figure}}    
\renewcommand{\thefigure}{{S\arabic{figure}}}%
\setcounter{equation}{0}
\renewcommand{\theequation}{S\arabic{equation}}

\section*{Supplementary Material}

\section{Proof of Theorem \ref{thm:explicit_discrete}}

\begin{theorem}
    For $w\in\mathbb R^{k+1}$ with pairwise distinct entries and $w_{\max} =
    \|w\|_{\infty}$, the Laplace transform of $\|S_{n,k}\|_{1,w}$ is given by
    \begin{multline}
        \mathbb E e^{t\|S_{n,k}\|_{1,w}} = k(-1)^k\cdot e^{t n w_{\max}} \times 
        \\
        \sum_{j=1}^{k} a_j^{e^{t(w-w_{\max})}} \bigg[ b_{n,k}
        \left( 1 - e^{t(w_j - w_{\max})(n+k-1)}\right)
        + \sum_{m=0}^{k-2} c_{n,k,m} \left( 1 - e^{t(w_j - w_{\max})}
        \right)^{k-1-m} \bigg],
        \label{eq:explicit_discrete_appendix}
    \end{multline}
    where for any $r\in\mathbb R^{k+1}$, $a_j^r = \prod_{m\neq j} (r_j-r_m)^{-1}$ and  
    \begin{align*}
        &b_{n,k} = \frac{(-1)^{k+1}}{n+k}\cdot\binom{n+k-1}{k-1}^{-1},
        &c_{n,k,m} = \frac{(-1)^m}{n+1}\cdot
        \frac{\binom{k-1}{m}}{\binom{n+m+1}{m}}.
    \end{align*}
\end{theorem}

\begin{proof}
    We observe that $S_{n,k} \sim \operatorname{Multinomial}
    (n,S_k)$, where $S_k\sim \operatorname{Dirichlet}(\mathbf{1}_{k+1})$, and so 
    \[
        \mathbb E e^{t\| S_{n,k} \|_{1,w}} = \mathbb E\mathbb E\left[ e^{t\| S_{n,k}
        \|_{1,w}} \mid S_k \right] = \mathbb E\left( \sum_{j=1}^{k+1} S_k\llbracket j
        \rrbracket e^{tw_j} \right)^n = \mathbb E \| S_k \|_{1,e^{tw}}^n,
    \]
    with $e^{tw}\in\mathbb R^{k+1}$ denoting the vector $(e^{tw_1}, ...,
    e^{tw_{k+1}})$. That is, the Laplace transform of interest is nothing but the
    $n^{\text{th}}$ moment of $\| S_k \|_{1,e^{tw}}$, which can be computed explicitly
    using the closed-form expression provided by Theorem \ref{thm:explicit_continuous}.
    This computation is lengthy, but straightforward, and results in
    \eqref{eq:explicit_discrete_appendix} as desired.
\end{proof}

\section{Proof of Theorem \ref{thm:explicit_continuous}}

\begin{theorem}
    As $n\to\infty$ with $k$ remaining fixed, $\mathbb P\left( \| S_{n,k}
    \|_{1,w} \leq x\right) = \mathbb P\left( \| S_k \|_{1,w} \leq x \right) +
    \varepsilon(x)$, where $S_k \sim
    \operatorname{Dirichlet}(\mathbf{1}_{k+1})$
    (with $\mathbf{1}_{k+1}\in\mathbb Z^{k+1}$ being the all-ones vector)
    is a uniform variable on the $k$-dimensional simplex, and $\|
    \varepsilon\|_{\infty}\in O(n^{-1})$. Moreover, with $a_j^w$ as in Theorem
    \ref{thm:explicit_discrete},
    \begin{equation}
        \mathbb P\left( \| S_k \|_{1,w} \leq x \right) = (-1)^{k}\sum_{j=1}^{k}
        a_{j}^w\left[ (x-w_j)_+
        \right]^k,
    \label{eq:explicit_continuous_appendix}
    \end{equation}
    as long as the components of $w$ are pairwise distinct.
\end{theorem}

\begin{proof}
    The $O(n^{-1})$ convergence is a consequence of a more general lemma.
    \begin{lemma*}
        Let $F_{n,k}^{p,w},F_k^{p,w}$ be the cumulative distribution functions of
        $\|S_{n,k} \|_{p,w}^p$ and $\|S_{k} \|_{p,w}^p$, respectively. Then
        \begin{equation}
            \| F_{n,k}^{p,w} - F_k^{p,w} \|_{\infty} = O(n^{-1}),
            \label{eq:convergence_rate}
        \end{equation}
        for every fixed $k \geq 2$.
    \end{lemma*}
    \begin{proof}[Proof of lemma]
        We approach the proof geometrically, showing that uniform samples from the
        discretized simplex converge to uniform samples from the continuous simplex as
        the discretization becomes finer. To do so, define the lattice $\Lambda =
        \mathbb Z^k \cap H$ where $H = \{x \in \mathbb R^k \ : \ \sum_{j=1}^k x_j =
        0\}$, and denote by 
        \begin{equation*}
            E^t = \{x \in \Delta^{k-1} : \ \| x \|_{p,w}^p \leq t\} = \{
            \|S_k\|_{p,w}^p \leq t\}
            \label{eq:level_set}
        \end{equation*}
        the $t$-level set of $F_k^{p,\bfw}$, we observe that since the fundamental
        domain of $\Lambda$ has diameter\break $\| (k-1, -1, \dots, -1)\|_2 =
        \sqrt{k(k-1)}$, the number $L(E^t,n)$ of lattice points in $nE^t$ is bounded by
        \begin{align*}
            \left(n-\sqrt{k(k-1)}\right)^{k-1}\mathrm{Vol}_{\Lambda}\left(
            E^t \right) &\leq
            \mathrm{d}\left( \Lambda \right) L(E^t,n) \notag \\
            &\leq \left(n+\sqrt{k(k-1)}\right)^{k-1}\mathrm{Vol}_{\Lambda}
            \left( E^t \right).
        \end{align*}
    Thus, in particular,
    \begin{align}
        \mu_{\discreteS_{n,k}}\left( nE^t \right) -
        \mu_{\Delta^{k-1}}\left(
        E^t \right) &= \dfrac{L(E^t,n)}{\binom{n+k-1}{k-1}} -
        (k-1)!\mathrm{Vol}_{\Lambda}\left( E^t \right) \notag \\ 
        &\leq (k-1)! \mathrm{Vol}_{\Lambda}\left( E^t \right) \left[
        \left( 1 +
        \dfrac{k}{n}\right)^{k-1} - 1 \right] \notag \\
        &\leq\sqrt{k}\sum_{j=1}^{k-1} \binom{k-1}{j} \left( \frac{k}{n}
        \right)^j,
        \label{eq:upper_bound}
    \end{align}
    where using $\mathrm{Vol}_{\Lambda}\left( E^1 \right) =
    \sqrt{k}/(k-1)!$ as
    an upper bound for $\mathrm{Vol}_{\Lambda}\left( E^t \right)$ turns
    \eqref{eq:upper_bound} independent of $t$. Similarly, a uniform
    lower bound
    is given by
    \begin{align}
        \mu_{\discreteS_{n,k}}\left( nE^t \right) -
        \mu_{\Delta^{k-1}}\left(
        E^t \right) &\geq (k-1)!\mathrm{Vol}_{\Lambda}\left( E^t \right)
        \left[
        \left( 1-\frac{2k}{n+k-1} \right)^{k-1} - 1 \right] \notag \\
        &\geq \sqrt{k} \sum_{j=1}^{k-1} \binom{k-1}{j}
        \left(\dfrac{-2k}{n+k-1}
        \right)^j.
        \label{eq:lower_bound}
    \end{align}
    Combining \eqref{eq:upper_bound} and \eqref{eq:lower_bound} gives
    \eqref{eq:convergence_rate} as desired.
    \end{proof}
    To arrive at \eqref{eq:explicit_continuous_appendix} then, write $\Gamma_{k+1}\| S_k
    \|_{1,w} = \sum_{j=1}^{k+1} w_j \mathcal E_j$,
    where $\Gamma_{k+1}\sim \operatorname{Gamma}(k+1,1)$, and
    $\mathcal E_j$ are $iid$ exponential variables of rate $1$, independent of
    $\Gamma_{k+1}$. This is a sum of independent variables, and thus admits
    factorization of its Laplace transform
    \[
        \mathbb E e^{t\Gamma_{k+1}\| S_k \|_{1,w}} = \prod_{j=1}^{k+1} \mathbb E 
        e^{tw_j \mathcal E_j} = \prod_{j=1}^{k+1} \frac{1}{1-tw_j}.
    \]
    On the other hand,
    \[
        \mathbb E e^{t\Gamma_{k+1} \| S_k \|_{1,w}} = \sum_{m=0}^{\infty}
        \frac{t^m}{m!} \mathbb E\Gamma_{k+1}^m \mathbb E \| S_k \|_{1,w}^m =
        \mathbb E \frac{1}{\left( 1 - t\| S_k \|_{1,w} \right)^{k+1}},
    \]
    rephrasing the task of identifying $\| S_k \|_{1,w}$'s distribution as inverting
    the Stieltjes-type transform
    \[
        \rho_k^{f}(z) = \mathbb E \frac{1}{\left( 1 - z\| S_k \|_{1,w} \right)^{k+1}}
        = \prod_{j=1}^{k+1} \frac{1}{1-zw_j},
    \]
    where $f$ denotes the density of $\| S_k \|_{1,w}$ (suppressing the dependence on
    $k$ and $w$ in $f$'s notation, as this will cause no ambiguity). To begin doing
    so, we observe that $f$ is a piece-wise polynomial of degree $k-1$ and knot points
    given by $w$ (as can be seen from the geometric interpretation of $\| S_k
    \|_{1,w}$), and thus has as $(k-1)^{\text{st}}$ derivative $\sum_{j=1}^{k} c_j
    \mathbbm{1}_{[w_j, w_j+1]}$ for some coefficients $c_j$. A $(k-1)$-fold
    integration by parts of $\rho_k^{f}(z)$ therefore yields
    \begin{multline*}
        \rho_k^{f}(z) = \frac{(-1)^{k-1}}{z^{k-1}k!} \rho_{2}^{f^{(k-1)}}(z) =
        \frac{(-1)^{k-1}}{z^{k-1}k!}\sum_{j=1}^{k} c_j \int_{w_j}^{w_{j+1}} 
        \frac{1}{(1-z x)^2} \ \mathrm{d}x 
        \\
        = \frac{(-1)^{k-1}}{z^{k-1}k!}\left( \frac{c_k}{z(1-zw_{k+1})} -
        \frac{c_1}{z(1-zw_1)} + \sum_{j=2}^k \frac{c_{j-1}-c_j}{z(1-zw_j)} \right),
    \end{multline*}
    which is meromorphic around the poles $1/w_j$, and so allows extraction of the 
    coefficients $c_j$ as 
    \[
        c_j = (-1)^{k-1}k! \sum_{m=1}^j \operatorname{Res}_{z=1/w_m} z^{k-1} 
        \rho_k^f(z) = (-1)^k k!\sum_{m=1}^k a_m^w.
    \]
    Using these coefficients to determine $f^{(k-1)}$, and integrating $k$ times gives
    \eqref{eq:explicit_continuous_appendix} as desired.
\end{proof}

\section{Proof of \thmref{thm:explicit_discrete_moments}}

\begin{theorem}
    Let $G(x,y) = \sum_{m=0}^{\infty} \mathrm{Li}_{-pm}(x)y^m/m!$, where
    $\mathrm{Li}_s(x) = \sum_{j=1}^{\infty} j^{-s} x^j$ is the 
    polylogarithm function. Denoting by $[x^ny^m]
    P(x,y)$ the $(n,m)^{\text{th}}$ coefficient of a power series $P$ in
    $x$ and $y$, we
    have
    \begin{equation}
        \mathbb E \left(\|\bfS_{n,k}\|_{p,\bfw}^{p}\right)^m =
        \dfrac{m!}{\binom{n+k-1}{k-1}} [x^ny^m] \prod_{i=1}^k G\left( x, w_i y
        \right).
        \label{eq:explicit_discrete_moments_supp}
    \end{equation}
    In particular, the first $m$ moments of $\|S_{i,j}\|_{p,\bfw}^p$ for
    $(i,j) \in
    \{0, \dots, n\} \times \{1, \dots, k\}$ can be computed in
    $O\left(nm\cdot  (\log nm) \cdot (\log k)\right)$ time.
\end{theorem}

\begin{proof}
    We first expand the left-hand side of \eqref{eq:explicit_discrete_moments_supp} to
    find
    \begin{align}
        \mathbb E\left( \|\bfS_{n,k}\|_{p,\bfw}^p \right)^m &=
        \sum_{\bfmath{\sample}\in\discreteS_{n,k}} 
        \mathbb P(\bfS_{n,k} = \bfmath{\sample})\left( \sum_{j=1}^k w_j
        \sample_j^p
        \right)^m \notag \\
        &= \binom{n+k-1}{k-1}^{-1}
        \sum_{\bfmath{\sample}\in\discreteS_{n,k}}
        \sum_{\bfmath{\eta}\in\discreteS_{m,k}} \binom{m}{\eta_1, \dots, \eta_k}
        \prod_{j=1}^k w_j^{\eta_j} \sample_j^{\eta_j p} \notag \\
        &= \frac{m!}{\binom{n+k-1}{k-1}}
        \underbrace{\sum_{\bfmath{\eta}\in\discreteS_{m,k}} \left(
            \sum_{\bfmath{\sample}\in\discreteS_{n,k}} \prod_{j=1}^k
            \frac{(w_j
        \sample_j^p)^{\eta_j}}{\eta_j!} \right)}_{A_{n,k,m,w}},
        \label{eq:multinomial_expansion}
    \end{align}
    so it remains to show that $A_{n,k,m,w} = [x^ny^m]\prod_{j=1}^k
    G(x,w_jy)$.
    By definition of $\mathrm{Li}_x(x)$, we have for every fixed
    $\bfmath{\eta}\in\discreteS_{m,k}$
    \begin{equation*}
        \sum_{\bfmath{\sample}\in\discreteS_{n,k}}\prod_{j=1}^k
        \frac{w_j^{\eta_j}
            \sample_j^{p\eta_j}}{\eta_j!} = [x^n] \prod_{j=1}^k
            \frac{\mathrm{Li}_{-p\eta_j}(x)}{\eta_j!}w_j^{\eta_j},
    \end{equation*}
    and so
    \begin{align*}
        A_{n,k,m,w} &= [x^n] \sum_{\bfmath{\eta}\in\discreteS_{m,k}}
        \prod_{j=1}^k
        \frac{\mathrm{Li}_{-p\eta_j}(x)}{\eta_j!}w_j^{\eta_j} \notag \\
        &= [x^n]\left\{[y^m] \prod_{j=1}^k \left( \sum_{i=0}^{\infty}
        \frac{\mathrm{Li}_{-pi}(x)}{i!} (w_jy)^i \right) \right\} \notag
        \\
        &= [x^ny^m] \prod_{j=1}^k G(x,w_jy),
    \end{align*}
    as desired. The $O\left( nm \cdot(\log nm)\cdot(\log k)\right)$
    runtime is now a direct consequence of 
    computing the Cauchy product of $k$ bivariate degree-$(n,m)$
    polynomials 
    using the Fast Fourier Transform.
\end{proof}

\section{Proof of \thmref{thm:explicit_continuous_moments}}

\begin{theorem} 
    Let $Q_p(x) = \sum_{m=0}^{\infty} (pm)!x^m/m!$.  Then,
    \begin{equation}
        \mathbb E\left( \| \bfS_k \|_{p,\bfw}^p \right)^m =
        \dfrac{(k-1)!m!}{(pm+k-1)!} [x^m] \prod_{j=1}^k Q_p(w_jx).
        \label{eq:explicit_continuous_moments_supp}
    \end{equation}
    In particular, the first $m$ moments of $\| S_j \|_{p,\bfw}^p$ for
    $j\in\{1, \dots, k\}$ can be computed in $O(m\cdot (\log m)\cdot
    (\log k))$ time.
\end{theorem}
\begin{proof}
    As in \eqref{eq:multinomial_expansion}, we expand the left-hand side
    of \eqref{eq:explicit_continuous_moments_supp} to obtain
    \begin{align}
        \mathbb E\left( \|\bfS_k\|_{p,\bfw}^p \right)^m &=
        \int_{\Delta^{k-1}}
        \left(\|\bfx\|_{p,\bfw}^{p}\right)^m \
        \mathrm{d}\mu_{\Delta^{k-1}}(\bfx)
        \notag \\
        &= \sum_{\bfmath{\eta}\in\discreteS_{m,k}}\binom{m}{\eta_1,
        \dots, \eta_k}
        \int_{\Delta^{k-1}} \prod_{j=1}^k \left( w_j^{\eta_j}
        x_j^{p\eta_j}
        \right) \ \mathrm{d}\mu_{\Delta^{k-1}}(x) \notag \\
        &= \dfrac{(k-1)!m!}{\sqrt{k}}
        \sum_{\bfmath{\eta}\in\discreteS_{m,k}}
        \Bigg(\prod_{j=1}^k \dfrac{w_j^{\eta_j}}{\eta_j!}\Bigg)
        \int_{\Delta^{k-1}}
        \prod_{i=1}^k x_i^{p\eta_i} \ \mathrm{d} \sigma(x) \notag \\
        &= \dfrac{(k-1)!m!}{\sqrt{k}}
        \sum_{\bfmath{\eta}\in\discreteS_{m,k}} 
        \Bigg(\prod_{j=1}^k\dfrac{w_j^{\eta_j}}{\eta_j!}\Bigg) \times
        \nonumber \\
          & \quad \int_{\Pi\Delta^{k-1}}
        \left(\prod_{i=1}^{k-1} x_i^{p\eta_i}\right)(1-x_1- \cdots -
        x_{k-1})^{p\eta_k} \sqrt{k} \ 
        \mathrm{d}\lambda_{k-1}(x) \label{eq:dirichlet_integral} \\
        &= \dfrac{(k-1)!m!}{(pm+k-1)!}
        \sum_{\bfmath{\eta}\in\discreteS_{m,k}}
        \prod_{j=1}^k \dfrac{(p\eta_j)!}{\eta_j!} w_j^{\eta_j}
        \label{eq:dirichlet_partition_function} \\
        &= \dfrac{(k-1)!m!}{(pm+k-1)!} [x^m]\prod_{j=1}^k \left(
        \sum_{i=0}^{\infty} \dfrac{(pi)!}{i!} (w_jx)^i\right)
        \label{eq:greenwood_multinomial}
    \end{align}
    where $\sigma(\mathrm{d}x)$ is (unnormalized) surface measure on
    $\Delta^{k-1}$, $\Pi\Delta^{k-1}$ the projection of $\Delta^{k-1}$
    on the hyperplane spanned by the first $k-1$ coordinate axes, 
    and \eqref{eq:dirichlet_partition_function} follows
    from recognizing the integral in \eqref{eq:dirichlet_integral} as the
    partition function of a Dirichlet variable with parameters $(p\eta_1, \dots,
    p\eta_k)$. We identify \eqref{eq:greenwood_multinomial} as
    \eqref{eq:explicit_continuous_moments_supp}, and thus complete the first part of
    the proof. The second part now follows as in Theorem
    \ref{thm:explicit_discrete_moments} from computing
    \eqref{eq:greenwood_multinomial} using the Fast Fourier Transform.
\end{proof}

\section{Proof of \propref{prop:hausdorff}}

\begin{proposition}
    Let $X\in[0,1]$ be either (a) continuous with density $f\in C^1\left([0,1]\right)$,
    or (b) discrete with support $\mathrm{supp}_X = \{x_0, \dots, x_N\}$. Moreover,
    denote by $B_{n,x}$ the degree-$n$ Bernstein polynomial approximating
    $\mathbbm{1}_{[0,x]}$. Then, for any resolution $\varepsilon_n \to 0, 
    \varepsilon_n > n^{-1/2}$, there exists $n_0(f,\varepsilon) \in 
    \mathbb N$, so that for all $n\geq n_0$,
    \begin{align}
        \sup_{x\in[0,1]} \big| \mathbb EB_{n,x}(X) - F(x) \big| &\leq
        \frac{\|f\|_{\infty} + 2\|f'\|_{\infty}+2}{n+1} \tag{a},
        \label{eq:bernstein_convergence_i_appendix} \\
        \sup_{x\in [0,1]\setminus \mathrm{supp}_X^{\varepsilon_n}} \big
        | \mathbb E B_{n,x}(X) - F(x) \big| &\leq e^{-2n\varepsilon_n^2}, 
        \tag{b} \label{eq:bernstein_convergence_ii_appendix}
    \end{align}
    where $\mathrm{supp}_X^{\varepsilon} = \left\{ x\in[0,1] ~:~ d(x, \mathrm{supp}_X)
    < \varepsilon \right\}$ is the $\varepsilon$-fattening of $\mathrm{supp}_X$.
\end{proposition}
\begin{proof}
    We first tackle \eqref{eq:bernstein_convergence_i_appendix} by recalling that the degree-$n$
    approximation by Bernstein polynomials \autocite{bernstein1912demo} is nothing but
    \begin{equation*}
        \mathbb EB_{n,x}(X) = \mathbb E\sum_{j=0}^{n-1}
        \mathbbm{1}_{\frac{j}{n} \leq x} \binom{n}{j} X^k\left( 1-X
        \right)^{n-k}.
    \end{equation*}
    To compute its approximation error, we choose a
    threshold $\varepsilon_n\to 0$ and investigate
    \begin{align}
        F(x) - \mathbb EB_{n,x}(X) &= \mathbb E\left(
        \mathbbm{1}_{[0,x]}(X) -
        B_{n,x}(X) \right) \notag \\
        &= \underbrace{\int_{[0,1]\setminus \{x\}^{\varepsilon_n}}
        \left(
        \mathbbm{1}_{[0,x]}(y) - B_{n,x}(y) \right)f(y) \
        \mathrm{d}y}_{A_{n,x}} 
        \notag \\
        &\quad + \underbrace{\int_{\{x\}^{\varepsilon_n}} \left(
        \mathbbm{1}_{[0,x]}(y) -
        B_{n,x}(y) \right) f(y) \ \mathrm{d}y}_{A_{n,x}'},
        \label{eq:integral_decomposition}
    \end{align}
    in which we treat the term $A_{n,x}$ first: Interpreting
    $B_{n,x}(y)$ as
    $\mathbb P\left( S_{n,y} \leq nx \right)$, where $S_{n,y} \sim
    \operatorname{Binomial}(n,y)$, we see that by standard large deviation
    estimates
    and Pinsker's inequality 
    \begin{align}
        \left|A_{n,x}\right| \leq (x&-\varepsilon_n) 
        \|f\|_{\infty}e^{-nD_{\text{KL}}\left( x \mid x-\varepsilon_n
        \right)} \notag \\ 
        &+ \|f\|_{\infty}(1-x+\varepsilon_n)e^{-n D_{\text{KL}} \left( x
        \mid x + 
        \varepsilon_n \right)} \leq
        \|f\|_{\infty}e^{-2n\varepsilon_n^2},
        \label{eq:pinsker}
    \end{align}
    where $D_{\text{KL}}\left( p \mid q \right)$ is the Kullback–Leibler
    divergence 
    (or the relative entropy) between
    a $\mathrm{Bernoulli}(p)$ and $\mathrm{Bernoulli}(q)$ distribution.
    To
    control $A_{n,x}'$ then, we Taylor expand $f$ to rewrite the
    integral in
    \eqref{eq:integral_decomposition} as
    \begin{align*}
        A_{n,x}' &= \int_{\{x\}^{\varepsilon_n}}\left(
        \mathbbm{1}_{[0,x]}(y) -
        B_{n,x}(y) \right)\left( f(x) + f'(\xi_{y,x})(y-x) \right) \
        \mathrm{d}y
        \notag \\
        &= \left(f(x) - M_n\cdot x\right)
        \underbrace{\int_{\{x\}^{\varepsilon_n}} \left(
        \mathbbm{1}_{[0,x]}(y) -
        B_{n,x}(y) \right) \mathrm{d} y}_{A_{n,x}''} \notag \\
        &\quad + M_n \underbrace{\int_{\{x\}^{\varepsilon_n}} \left(
        \mathbbm{1}_{[0,x]}(y) - B_{n,x}(y) \right)y \
        \mathrm{d}y}_{A_{n,x}'''},
    \end{align*}
    where $\min_{y\in\{x\}^{\varepsilon_n}} f'(y) \leq M_n \leq
    \max_{y\in\{x\}^{\varepsilon_n}} f'(y)$. In particular, since we
    assumed
    $f\in C^1\left( [0,1] \right)$ and $\varepsilon_n\to 0$, there must
    exist a $n_0'$ 
    so that $f'(x)-1 \leq M_n \leq f'(x) + 1$ for all $n\geq n_0'$. So
    it remains
    to control $A_{n,x}''$ and $A_{n,x}'''$, which can be done in a
    manner
    similar to \eqref{eq:pinsker}:
    
    \begin{align}
        \left|A_{n,x}''\right| &\leq \int_{[0,1]} \left(
        \mathbbm{1}_{[0,x]}(y) -
        B_{n,x}(y) \right) \ \mathrm{d}y + e^{-2n\varepsilon_n^2} =
        \frac{x-1}{n+1} + e^{-2n\varepsilon_n^2}
        \notag  \\
        &\leq \frac{1}{n+1} + e^{-2n\varepsilon_n^2} \notag \\
        \left| A_{n,x}'''\right| &\leq \int_{[0,1]} \left(
        \mathbbm{1}_{[0,x]}(y) - B_{n,x}(y) \right) y \ \mathrm{d}y +
        e^{-2n\varepsilon_n^2} \label{eq:a_triple_prime} \\
        &= \frac{3nt(x-1) + 2(x^2-1)}{2(n+1)(n+2)} +
        e^{-2n\varepsilon_n^2} \leq
        \frac{1}{n+1} + e^{-2n\varepsilon_n^2} \notag,
    \end{align}
    provided $n\geq 4$. Finally, combining
    \eqref{eq:integral_decomposition}-\eqref{eq:a_triple_prime}, we
    obtain
    \begin{equation*}
        \left| \hat{F}_n(x) - F(x) \right| \leq \frac{\|f\|_{\infty} +
        2\|f'\|_{\infty} + 2}{n+1} + 2\left( \|f\|_{\infty} +
        \|f'\|_{\infty} \right)e^{-2n\varepsilon_n^2},
    \end{equation*}
    independently of $x$. Choosing $\varepsilon_n \geq
    n^{-\frac{1}{2}+\delta}$ and 
    $n_0$ so large that the first term dominates
    the second yields \eqref{eq:bernstein_convergence_i_appendix}. (ii) follows
    in a very
    similar manner by observing that for $n$ such that $\varepsilon_n <
    h$,
    any $x\in[0,1]\setminus \mathrm{supp}_X^{\varepsilon_n}$ satisfies
    \begin{equation*}
        \left|\mathbbm{1}_{[0,x]}(y) - B_{n,x}(y)\right| \leq
            e^{-2n\varepsilon_n^2}.
    \end{equation*}
    Therefore,
    \begin{align*}
        \left|F(x) - \mathbb EB_{n,x}(X) \right| &\leq
        \sum_{y\in\mathrm{supp}_X} 
        \mathbb P\left( X = y \right) \left| \mathbbm{1}_{[0,x]}(y) -
        B_{n,x}(y)
        \right| \notag \\
        &\leq e^{-2n\varepsilon_n^2}
    \end{align*}
    which is \eqref{eq:bernstein_convergence_ii_appendix}.
\end{proof}

\section{Proof of Proposition \ref{prop:right_tail}}

\begin{proposition}
    Without loss of generality, assume $w\in\mathbb R_+^{k+1}$ and $\|w\|_{\infty} =
    1$, and denote by
    \[
        W_{w} = \left| \{ 1 \leq j \leq k+1 ~ : ~ w_j = 1 \} \right|
    \]
    the
    number of weight components assuming value $1$. Then the density $f_k^{p,w}$ of 
    $\| S_k \|_{p,w}^p$ is analytic on $[x_0,1]$, where
    \[
        x_0 = 
        \begin{cases}
            \frac{1}{2^{p-1}} &\text{ if }W_w = k+1 \\
            \max_{j: w_j < 1} w_j &\text{ otherwise},
        \end{cases}
    \]
    and its degree $r$ Taylor polynomial around $1$ can be computed in $O\left(
    \frac{r}{p}\log\frac{r}{p} \log k + [r\log r]^2\right)$ time. For $r=k-2$ it reads
    \begin{equation*}
        f_k^{p,w}(x) = \dfrac{(k-1)W_{w}}{2^{k-1}}\left( 1-x \right)^{k-2} +
        O\left( (1-x)^{k-1} \right).
        \label{eq:general_taylor_appendix}
    \end{equation*}
    In particular, Greenwood's statistic satisfies
    \begin{equation*}
        f_k^{2,\mathbf{1}_{k}}(x) = \dfrac{\binom{k}{2}}{2^{k-2}} 
        \left( 1-x \right)^{k-2} + O\left( (1-x)^{k-1} \right).
        \label{eq:taylor_expansion_appendix}
    \end{equation*}
\end{proposition}

\begin{proof}
    $f_{k}^{p,w}$ being analytic around $[x_0,1]$ follows directly from the geometric 
    perspective that has been used extensively in previous proofs already. 
    The asymptotic behavior of its moments governs $f_k^{p,w}$ on this interval. The
    following result clarifies this behavior.
    \begin{lemma*}
        For $p\geq 2$ and $k\geq 2$, and fixed weights $w_i\in[0,1]$, for
        all $i\in [k]$, we have
        \begin{equation}
            \lim_{m\to\infty} m^{k-1} \left(\mathbb
            E\|\bfS_k\|_{p,\bfw}^p\right)^m =
            \dfrac{(k-1)!}{p^{k-1}} \cdot W_{\bfw},
            \label{eq:moment_rates}
        \end{equation}
        where $W_{\bfw} = \left| \{ 1 \leq i \leq k ~ : ~ w_i = 1 \}
        \right|$ is the number of weights taking value $1$. In particular, the Greenwood
        statistic satisfies
        \begin{equation*}
            \lim_{m\to\infty} m^{k-1} \left( \mathbb E
            \|\bfS_k\|_{2,\bfmath{1}_k}^2
            \right)^m = \dfrac{k!}{2^{k-1}}.
        \end{equation*}
    \end{lemma*}
    \begin{proof}[Proof of lemma]
        We first rewrite \eqref{eq:dirichlet_partition_function} as
        \begin{equation}
            \mathbb E\left(\|\bfS_k\|_{p,\bfw}^p \right)^m =
            \dfrac{1}{\binom{pm+k-1}{k-1}}
            \sum_{\bfmath{\eta}\in\discreteS_{m,k}}
            \dfrac{\binom{m}{\eta_1, \dots, \eta_k}}{ \binom{pm}{p\eta_1,
            \dots,
            p\eta_k}}\prod_{j=1}^kw_j^{\eta_j} =
            \dfrac{1}{\binom{pm+k-1}{k-1}} s_m^{\bfw},
            \label{eq:rewritten_decay}
        \end{equation}
        which has leading order $O\left( m^{-(k-1)} \right)$, if we can show that
        $s_m^{\bfw}$ is $\Omega(1)$. To do so, we proceed by induction on
        $k$, the
        length of $w$, proving  that in fact $\lim_{m\to\infty} s_m^{\bfw} =
        W_{\bfw}$. It is
        straightforward to check that for $\eta \in \{2,\dots,m-1\}$,
        $\binom{m}{\eta} /
        \binom{pm}{p\eta}$ is bounded above by $\binom{m}{2} /
        \binom{2m}{2\eta}$,
        and thus for the base case $k=2$ we have
        \begin{align}
            s_m^{(w_1,w_2)} = \sum_{\eta=0}^m
            \dfrac{\binom{m}{\eta}}{\binom{pm}{p\eta}}
            w_1^{\eta}w_2^{m-\eta} \leq w_1^m + w_2^m &+
            \dfrac{\binom{m}{1}}{\binom{pm}{p}} +
            (m-2)\dfrac{\binom{m}{2}}{\binom{pm}{2p}} \notag \\
            &\xrightarrow{m\to\infty} \mathbbm{1}_{w_1 = 1} +
            \mathbbm{1}_{w_2 = 1}
            = W_{(w_1,w_2)},
            \label{eq:base_case}
        \end{align}
        as desired. For the inductive step, we condition on the first entry
        of
        $\eta$ to obtain
        \begin{align*}
            s_m^{(w_1, \dots, w_k)} &= \sum_{\ell = 0}^m 
            \dfrac{\binom{m}{\ell}}{\binom{pm}{p\ell}} 
            w_1^{\ell} \sum_{\bfmath{\eta}\in\discreteS_{m-\ell,k-1}} 
            \dfrac{\binom{m-\ell}{\eta_1, \dots, \eta_{k-1}}}{ 
            \binom{p(m-\ell)}{p\eta_1,\dots,p\eta_{k-1}}} \prod_{j=1}^{k-1} 
            w_{j+1}^{\eta_j} \notag \\
            &= s_m^{(w_2, \dots, w_k)} + w_1^m + O\left( m^{-1} \right)
            \notag \\
            &\xrightarrow{m\to\infty} W_{(w_2, \dots, w_k)} +
            \mathbbm{1}_{w_1 = 1}
            = W_{(w_1, \dots, w_k)},
        \end{align*}
        where we used the inductive hypothesis on $s_m^{(w_1, \dots, w_k)}$,
        and as
        in \eqref{eq:base_case}, bounded summands corresponding to $\ell \in
        \{2,
        \dots, m-1\}$ by $\binom{m}{2}/\binom{pm}{2p}$.
        The lemma now follow
        from taking the limit as $m\to\infty$ in \eqref{eq:rewritten_decay}.
    \end{proof}
    Let $f^{p,\bfw}_k(x) = \sum_{j=0}^{\infty} c_j (1-x)^j$ be the
    Taylor expansion
    of $f^{p,\bfw}_k$ around $x_0 = 1$. We first notice that for any $r
    \geq 0$,
    \begin{equation*}
        \int_0^1 x^m(1-x)^r ~ \mathrm{d}x = \dfrac{1}{m+r+1}\cdot
        \dfrac{1}{\binom{m+r}{r}},
    \end{equation*}
    and hence, using the fact that $f_k^{p,\bfw}$ is bounded,
    \begin{align*}
        \mathbb E\left( \|\bfS_{n,k}\|_{p,\bfw}^p \right)^m &=
        \int_{0}^1 x^m
        f^{p,\bfw}_k(x) ~ \mathrm{d}x + O\left(e^{-m}\right) \notag \\
        &= \sum_{j=0}^{\infty} c_j \int_0^1 x^m(1-x)^j ~ \mathrm{d}x + 
        O\left( e^{-m} \right)\notag \\ &
        = \sum_{j=0}^{\infty} c_j \dfrac{1}{m+j+1}
        \dfrac{1}{\binom{m+j}{j}} +
        O\left( e^{-m} \right).
    \end{align*}
    Identifying the $(k-2)^{\text{nd}}$ term with \eqref{eq:moment_rates}
    immediately yields the first-order Taylor expansions of $f_k^{p,w}$.
    
    To compute higher-order expansion, we recall from \eqref{eq:rewritten_decay} that
    $\mu_m = \mathbb E \left( \|\bfS_k\|_{p,\bfw} \right)^{pm}$ can be written as
    \begin{equation}
        \mu_m = \dfrac{1}{\binom{pm+k-1}{k-1}} \sum_{\bfmath{\eta}\in\discreteS_{m,k}} \dfrac{\binom{m}{\eta_1, \dots, \eta_k}}{\binom{pm}{p\eta_1, \dots, 
        p\eta_k}} \prod_{j=1}^k w_j^{\eta_j} = \dfrac{s^{\bfw}_m}{\binom{pm+k-1}{k-1}},
        \label{eq:mu_LHS}
    \end{equation}
    where $s_m^{\bfw} = \sum_{j=0}^{\infty} \sigma_j^{\bfw}(m)\cdot m^{-j}$ with
    $\sigma_j^{\bfw}(m)$ remaining constant $\sigma_j^{\bfw}$ past some threshold
    $m_j^w$. We also have
    \begin{align}
        \mu_m = \int_0^1 x^mf_k^{p,\bfw}(x) \ \mathrm{d}x &=
        \sum_{j=0}^{\infty} c_j^{\bfw}
        \int_0^1 x^m(1-x)^j \ \mathrm{d}x + O\left( e^{-m} \right) \notag \\ 
        &= \sum_{j=0}^{\infty} c_j^{\bfw} \left[ (m+j+1) \binom{m+j}{j} \right]^{-1}
        + O\left( e^{-m} \right),
        \label{eq:mu_RHS}
    \end{align}
    which suggests that by matching coefficients in \eqref{eq:mu_LHS} and
    \eqref{eq:mu_RHS} we should be able to translate between
    $\sigma_j^{\bfw}$ and $c_j^{\bfw}$. For this to be helpful, we need to understand
    $\sigma_j^{\bfw}$:
    \begin{lemma*}[$\sigma_j^{\bfw}$ recursion] 
        Defining $b^{j}_r = \left[m^{-r}\right] \binom{m}{j} /
        \binom{pm}{pj}$ and employing notation as in
        \eqref{eq:mu_LHS}, we have 
        \begin{equation}
            \sigma_r^{\bfw} = \sum_{j=0}^{r'} \left( \sigma_j^{(w_k,0)} \cdot
            \sigma^{\bfw_{-k}}_{r-j} + \mathbbm{1}_{w_k=1} s^{\bfw_{-k}}_j
            \cdot b^j_r \right),
            \label{eq:sigma_coefficients}
        \end{equation}
        with initial condition $\sigma_r^{w_1,w_2} = \sum_{j=0}^{r'}
        b^j_r\left( \mathbbm{1}_{w_2=1 }w_1^j + \mathbbm{1}_{w_1=1} w_2^j\right)$,
        where $r' = \lfloor r/(p-1) \rfloor$. In particular, we can compute
        $\sigma_{r}^{\bfw}$ in $O\left( r'\log{r'}\log{k} +
        \left[r\log{r}\right]^2 \right)$ time.
    \end{lemma*}
    \begin{proof}[Proof of lemma]
        Slightly abusing notation, we have
        \begin{align*}
            \sigma_r^{\bfw} &= \left[ m^{-r} \right] s_m^{\bfw} = \left[ m^{-r}
                \right] \sum_{\bfmath{\eta}\in\discreteS_{m,k}}
                \dfrac{\binom{m}{\eta_1, \dots,
                \eta_k}}{\binom{pm}{p\eta_1, \dots, p\eta_k}} \prod_{j=1}^k
                w_j^{\eta_j} \notag \\
                &= \left[ m^{-r} \right]\sum_{\omega=0}^m
                \dfrac{\binom{m}{\omega}}{\binom{pm}{p\omega}} w_k^{\omega}
                \sum_{\bfmath{\eta}\in\discreteS_{m-\omega,k-1}}
                \dfrac{\binom{m-\omega}{\eta_1, \dots,
                \eta_{k-1}}}{\binom{p(m-\omega)}{p\eta_1, \dots,
                p\eta_{k-1}}} \prod_{j=1}^{k-1}w_j^{\eta_j}\notag\\
                & = \left[ m^{-r}
                \right] \sum_{\omega=0}^{m} \dfrac{\binom{m}{\omega}}{ 
                \binom{pm}{p\omega}} w_k^{\omega} \cdot
                s_{m-\omega}^{\bfw_{-k}} \notag \\
                &= \left[ m^{-r} \right]\sum_{\omega=0}^{r'}
                \dfrac{\binom{m}{\omega}}{\binom{pm}{p\omega}} w_k^{\omega}
                \cdot s_{m}^{\bfw_{-k}} + \left[ m^{-r} \right]
                \sum_{\omega=0}^{r'} \dfrac{\binom{m}{\omega}}{\binom{
                pm}{\omega}} w_k^{m-\omega} \cdot
                s_{\omega}^{\bfw_{-k}} \notag \\
                &= \left[ m^{-r} \right] \sum_{\omega=0}^{r'} s_m^{\bfw_{-k}}
                \sum_{j=0}^{\infty} b^{\omega}_j w_k^{\omega} m^{-j} 
                + \left[ m^{-r} \right] \sum_{\omega=0}^{r'} w_k^{m-\omega}
                s_{\omega}^{\bfw_{-k}}  \sum_{j=0}^{\infty} b^{\omega}_j
                m^{-j}  \notag \\
                &= \sum_{j=0}^{r'} \sigma_{r-j}^{\bfw_{-k}}\cdot
                \sigma_j^{w_k,0} + \mathbbm{1}_{w_k = 1} \sum_{\omega=0}^{r'} 
                s_{\omega}^{w_{[1:k-1]}}\cdot b_j^{\omega} \\
                & =
                \sum_{j=0}^{r'}\left( \sigma_j^{w_k,0}\cdot
                \sigma_{r-j}^{\bfw_{-k}} + \mathbbm{1}_{w_k = 1}
                s_j^{\bfw_{-k}}\cdot b^j_r \right),
        \end{align*}
        as desired. To see that \eqref{eq:sigma_coefficients} can be
        computed in $O\left( r'\log{r'}\log{k} + \left[ r\log{r} \right]^2
        \right)$ time, we notice that calculation of $s_r^{\bfw_{-k}}$ is
        $O\left( r'\log{r'}\log{k} \right)$ by the same reasoning as in
        \textsc{Proposition} 2, and $b^{j}_r$, written as,
        \begin{align*}
            b^j_r &= \left[ m^{-r} \right]
            \dfrac{\binom{m}{j}}{\binom{pm}{pj}} = (pj-1)!_{p} \left[
            m^{-r} \right] \prod_{\substack{\ell = p(m-j)+1 \\ p ~ \nmid ~ 
            \ell}}^{pm-1} \dfrac{1}{pm}\cdot \dfrac{1}{1-\frac{\ell}{pm}}
            \notag \\ 
            &= (pj-1)!_{p}\left[ m^{-r} \right] \prod_{\substack{\ell =
            p(m-j)+1 \\ p ~ \nmid ~ \ell}}^{pm-1} R\left(
            \frac{\ell}{pm} \right),
        \end{align*}
        where $R(x) = \sum_{j=0}^{\infty} x^j$ is again a convolution of
        $(p-1)\cdot r' = r$ polynomials and hence computable in 
        $O\left( \left[ r\log{r} \right]^2 \right)$.
    \end{proof}

    With a proper understanding of $\sigma_j^{\bfw}$ at hand, we may rewrite \eqref{eq:mu_LHS} as
    \begin{equation}
        \mu_m = \sum_{j=0}^{\infty} \left( \sum_{\omega=0}^j
        a^k_{\omega}\cdot \sigma_{j-\omega}^{\bfw} \right) m^{-j} + O\left(
        e^{-m} \right),
        \label{eq:mu_LHS_coefficients}
    \end{equation}
    where $a^k_{\omega} = \left[ m^{-\omega} \right]
    \binom{pm+k-1}{k-1}^{-1}$. Similarly, expanding \eqref{eq:mu_RHS}
    yields
    \begin{equation}
        \mu_m = \sum_{j=0}^{\infty}\left( \sum_{\omega=0}^{j-1} d^{\omega}_j
        \cdot c^{\bfw}_{\omega} \right) m^{-j} + O\left( e^{-m} \right),
        \label{eq:mu_RHS_coefficients}
    \end{equation}
    where $d^{\omega}_j = \left[ m^{-j} \right] \left[ (m+\omega+1)
    \binom{m+\omega}{\omega} \right]^{-1}$. Consequently, matching the
    $r^{\text{th}}$ coefficients in \eqref{eq:mu_LHS_coefficients} and
    \eqref{eq:mu_RHS_coefficients} allows to solve for $c^{\bfw}_{r}$:
    \begin{equation*}
        c^{\bfw}_r = \dfrac{1}{r!}\left[ \sum_{j=k-1}^{r+1} a^k_j\cdot
        \sigma^{\bfw}_{r+1-j} - \sum_{j=k-2}^{r-1} d^j_{r+1} c^{\bfw}_j \right],
    \end{equation*}
    where in the choice of summation indices we have used the fact that
    $a^k_j = 0$ for $j\in\{0, \dots, k-2\}$ and $c^{\bfw}_j = 0$ for $j\in\{0,
    \dots k-3\}$. Now, $\{d^0_{r+1}, \dots,
    d^{r-1}_{r+1}\}$ can be found in $O\left(r\left( \log{r}
    \right)^2\right)$ time, and given $a,b$ and $d$, the recursion is
    solved in $O\left(r^2\right)$ steps, amounting to a total complexity of\break
    $O\left( r\left( \log{r} \right)^2 + r^2 + r'\log{r'}\log{k} + \left[
    r\log{r} \right]^2 \right)
    = O\left( r'\log{r'}\log{k} + \left[
    r\log{r} \right]^2 \right)$.
\end{proof}

\section{Proof of Proposition \ref{prop:multivariate_extension}}

\begin{proposition}
    For $m$ weights $w^1, ..., w^m\in\mathbb R_{k+1}$, each of pairwise distinct entries,
    the Laplace transform of the tuple $S_{(m)} = \left( \|S_{n,k}\|_{1,w^1}, ..., \|
    S_{n,k} \|_{1,w^m} \right)$ is given by
    \begin{multline*}
        \mathbb E e^{\langle t, S_{(m)} \rangle} = k(-1)^k\cdot e^{n w^{\max}} \times 
        \\ \sum_{j=1}^{k} a_j^{e^{\omega_j} } \bigg[ b_{n,k} \left( 1 - e^{\omega_j
        (n+k-1)} \right) + \sum_{m=0}^{k-2} c_{n,k,m} \left( 1 - e^{\omega_j}
        \right)^{k-1-m} \bigg],
        \label{eq:multivariate_extension}
    \end{multline*}
    where $\omega_j = \sum_{r=1}^m t_r w_j^r - w^{\max}$ with $w^{\max} = \max_j
    \sum_{r=1}^m t_r w^r_j$, and $a_j^w, b_{n,k}$ and $c_{n,k,m}$ as in Theorem
    \ref{thm:explicit_discrete}.
    
    Moreover, the joint moments of $S_{(m)}$ can be computed in $O\big(n\prod_{j=1}^r
    m_j$ $\times  (\log n\prod_{j=1}^r m_j) \times k\big)$ time as
    \begin{equation*}
        \mathbb E\prod_{j=1}^r \left(\|\bfS_{n,k}\|_{p_j,\bfw^j}^{p_j}\right)^{m_j} =
        \dfrac{\prod_{j=1}^r m_j!}{\binom{n+k-1}{k-1}} [x^ny_1^{m_1}\cdots y_r^{m_r}]
        \prod_{i=1}^k
        G_r\left( x, w^1_i y_1, \dots, w^r_i y_r \right),
        \label{eq:main_mulitvariate_appendix}
    \end{equation*}
    where $G_r(x,y_1, \dots, y_r) = \sum_{m_1, \dots, m_r=0}^{\infty}
    \mathrm{Li}_{-\sum_{j=1}^r p_jm_j}(x)\prod_{j=1}^r y_j^m/m_j!$. These joint
    moments
    can be used to approximate $\mathbb P\left( \| S_{(m)} \|_{\infty} \leq x \right)$
    up to $\varepsilon$ accuracy in $O\left( \varepsilon^{-1} \right)$ time.
\end{proposition}

\begin{proof}
    The first part of the statement follows from the fact that 
    \[
        \langle t, S_{(m)} \rangle = \sum_{j=1}^{m} t_j \| S_{n,k} \|_{1,w^j} = \|
        S_{n,k} \|_{1, \sum_{j=1}^m t_jw^j},
    \]
    and following the same reasoning as in Theorem \ref{thm:explicit_discrete}.
    Similarly, the second part closely follows the arguments of Theorem
    \ref{thm:explicit_discrete_moments}.
    \begin{align}
        \mathbb E\prod_{i=1}^r \left( \|\bfS_{n,k}\|_{p_i,\bfw_i}^{p_i}
        \right)^{m_i} 
        &= \sum_{\bfmath{\sample}\in\discreteS_{n,k}} 
        \mathbb P(\bfS_{n,k} = \bfmath{\sample})\prod_{i=1}^r \left(
        \sum_{j=1}^k w_{i,j} \sample_j^p
        \right)^{m_i} \notag \\
        &= \binom{n+k-1}{k-1}^{-1}
        \sum_{\bfmath{\sample}\in\discreteS_{n,k}}
        \prod_{i=1}^r \left[ \sum_{\bfmath{\eta}_i
        \in\discreteS_{m_i,k}} \binom{m_i}{\eta_{i,1}, \dots,
        \eta_{i,k}}
        \prod_{j=1}^k w_{i,j}^{\eta_{i,j}} \sample_j^{\eta_{i,j} p_i}
        \right] \notag \\
        &= \frac{\prod_{i=1}^r m_i!}{\binom{n+k-1}{k-1}} 
        \underbrace{
            \sum_{\bfmath{\eta_1}\in\discreteS_{m_1,k}}\cdots
            \sum_{\bfmath{\eta_r}\in\discreteS_{m_r,k}} 
            \left(
                \sum_{\bfmath{\sample}\in\discreteS_{n,k}} \prod_{i=1}^r
                \prod_{j=1}^k
            \frac{(w_{i,j} \sample_j^{p_i})^{\eta_{i,j}}}{\eta_{i,j}!}
            \right)
        }_{A_{n,k,m_i,w_i}},
        \label{eq:multivariate_multinomial_expansion}
    \end{align}
    so it remains to show that $A_{n,k,m_i,w_i} = [x^ny_1^{m_1}\cdots
    y_r^{m_r}]\prod_{j=1}^k
    G_r(x,w_{1,j}y_1, \dots, w_{r,j}y_r)$.
    By definition of $\mathrm{Li}_x(x)$, we have for every fixed
    $\bfmath{\eta}\in\prod_{i=1}^r\discreteS_{m_i,k}$
    \begin{equation}
        \sum_{\bfmath{\sample}\in\discreteS_{n,k}}\prod_{j=1}^k
        \frac{w_j^{\sum_i\eta_{i,j}}
        \sample_j^{\sum_ip_i\eta_j}}{\prod_i\eta_{i,j}!} = [x^n]
        \prod_{j=1}^k
        \frac{\mathrm{Li}_{-\sum_ip_i\eta_{i,j}}(x)}{\prod_i\eta_{i,j}!
        }w_j^{\sum_i\eta_{i,j}},
        \label{eq:multivariate_x_convolution}
    \end{equation}
    and so
    \begin{align}
        A_{n,k,m_i,w_i} &= [x^n]
        \sum_{\bfmath{\eta_1}\in\discreteS_{m_1,k}} \cdots 
        \sum_{\bfmath{\eta_r}\in\discreteS_{m_r,k}}\prod_{j=1}^k
        \frac{\mathrm{Li}_{-\sum_ip_i\eta_{i,j}}(x)}{\prod_i
        \eta_{i,j}!}w_j^{\sum_i\eta_{i,j}} \notag \\
        &= [x^n]\left\{[y_1^{m_1}]\cdots [y_r^{m_r}] \prod_{j=1}^k
        \left( \sum_{m_1, \dots, m_r=0}^{\infty}
        \frac{\mathrm{Li}_{-\sum_i p_i m_i}(x)}{\prod_i m_i!} \prod_i
        (w_{i,j} y_i)^{m_i} \right) \right\} \notag \\
        &= [x^ny_1^{m_1}\cdots y_r^{m_r}] \prod_{j=1}^k
        G_r(x,w_{1,j}y_1, \dots, w_{r,j}y_r),
        \label{eq:multivariate_y_convolution}
    \end{align}
    as desired. 
\end{proof}

\section{Proof of Proposition \ref{prop:klotz}}

\begin{proposition}
    For $a,b>0$, define random variable $X_{a,b}$ through their densities
    \[
        f_{a,b}(x) = 
        \left(
            \begin{cases}
                 \frac{1}{a} \frac{x+b}{b-a}    & \text{if }-b \leq x < -a \\
                 -\frac{1}{x}                   & \text{if }-a \leq x < -1 \\
                 1                              & \text{if }-1 \leq x < 1  \\
                 \frac{1}{x}                    & \text{if }1  \leq x < a  \\
                 \frac{1}{a} \frac{x-b}{a-b}    & \text{if }a  \leq x < b
            \end{cases}
        \right)/Z,
    \]
    where $Z = 2\log a + 1 + \kappa$, with $\kappa = b/a$. Then as $a\to\infty$ while
    keeping $\kappa\in o(\log a)$, $e_w^{\sigma}(F_{a,b})\to 0$ for any $w\in
    \mathcal C^1([0,1])$ whose derivative is bounded by $C ( x^{-1}|\log x|^p +
    (1-x)^{-1}|\log(1-x)|^p)$ for some constants $C>0$ and $p>0$.
\end{proposition}

\begin{proof}
    Assuming without loss of generality that $\int w = 0$ and $\int w^2 = 1$, it follows
    from \autocite{holst1980asymptotic} that under local scale alternatives $Y =
    (1+\sigma_n)X$,
    $\| S_{n,k} \|_{1,w}$, suitably standardized, is distributed $\mathcal N(\mu, 1)$,
    where
    \[
        \mu_S = \frac{1}{\sqrt{(1+\alpha) \operatorname{Var} w(U)}} \int_0^1
        w'(x) F^{-1}(x)f\left(F^{-1}(x)\right) \ \mathrm{d}x,
    \]
    with $f$ and $F$ the density and CDF of $X$, respectively, $U$ a uniform variable
    on $[0,1]$, and as long as $w$ is in $\mathcal C^1([0,1])$ and satisfies the
    boundary assumption given in the Proposition statement. Similarly, standard
    CLT-type computations for an appropriately normalized $F$-statistic $F_{n,k}$ show
    that it behaves asymptotically normal of unit variance and expectation
    \[
        \mu_F = \frac{2}{\sqrt{\left(\frac{\mu_4}{\sigma^4}-1\right)(1+\alpha)}},
    \]
    where $\mu_4$ and $\sigma^2$ are the fourth moment and variance of $X$,
    respectively. The Pitman efficiency between $\| S_{n,k} \|_{1,w}$ and $F_{n,k}$ is
    thus given by 
    \[
        e_w^{\sigma}(F) = \frac{\left(\frac{\mu_4}{\sigma^4}-1\right)}{2^2} \left(
        \int_0^1 w'(x) F^{-1}(x)f\left(F^{-1}(x)\right) \ \mathrm{d}x \right)^2.
    \]
    The task then is to show that this quantity can be made arbitrarily small for the
    type of $w$ in question and $F=F_{a,b}$. To do so, we first observe that
    the factor $(\mu_4/\sigma^4 - 1)/2^2$ is straightforwardly computed to be of order
    $O(\kappa^{-1}\log a)$, and so it suffice to show that
    \[
        \int_0^1 w'(x) F_{a,b}^{-1}(x)f_{a,b}\left(F_{a,b}^{-1}(x)\right) \ \mathrm{d}x
    \]
    is $o\left( \sqrt{\kappa/\log a} \right)$. We will demonstrate that this rate can
    be achieved for the $\int_{1/2}^1$ segment of this integral, from which the
    $\int_0^{1/2}$ component will follow by symmetry of $f_{a,b}$. The explicit form of
    $f_{a,b}$ allows computation of $F_{a,b}^{-1}\cdot \left( f_{a,b} \circ 
    F_{a,b}^{-1}\right)$ as
    \[
        F_{a,b}^{-1}(x)f_{a,b}\left(F_{a,b}^{-1}(x)\right) \quad
        \begin{cases}
             \leq 1/Z       &\text{if }x \leq x_{\kappa} \\
             = \sqrt{\frac{2(1-x)}{Z(\kappa-1)}} \left( \kappa -
             \sqrt{2Z(1-x)(\kappa-1)}
             \right)        &\text{if }x > x_{\kappa},
        \end{cases}
    \]
    where $x_{\kappa} = 1-(\kappa-1)/(2Z)$, and so bounding $w'$ on $[1/2,1)$ by
    $C(1-x)^{-1}\log^p(1-x)$, the integral of interest is bounded by the magnitude of
    \[
        \frac{x_{\kappa}}{Z} + 2C \left( \Gamma(p+1, \log\frac{2Z}{\kappa-1}) - 
        \frac{\kappa 2^{p+1} \Gamma(p+1,\frac{1}{2}\log\frac{2Z}{\kappa-1})}{
        \sqrt{Z(\kappa-1)}} \right) \in O\left(\frac{\kappa-1}{Z}\log
        \frac{Z}{\kappa-1}\right),
    \]
    which is in $o(\sqrt{\kappa/\log a})$ as long as $\kappa/Z \to 0$, which in turn
    is achieved whenever $\kappa\in o(\log a)$.
\end{proof}

\section{Examples}

As already indicated in the proof of Proposition \ref{prop:klotz}, it follows from
arguments of \autocite{holst1980asymptotic} (or general LAN theory as developed in [REF])
that $\| S_{n,k} \|_{1,w}^1$, under local alternatives $Y_n\sim G_n = F + \widetilde
H_n/\sqrt{n}$ and suitably normalized, behaves like a Gaussian $\mathcal N(\mu_h, 1)$
with expectation
\[
    \mu_h = \int_{0}^{1} w(x) h(x) \ \mathrm{d}x,
\]
where $h = \lim_{n\to\infty}\sqrt{n}(g^F_n - 1)$ with $g^F_n$ the density of $F(Y)$. 
Consistency and asymptotic power results can thus be read off from the magnitude of 
$\mu_h$; and doing so for, e.g., $\| S_{n,k} \|_{1,(0,1,...,k)/k}$ (corresponding to
the widely used Mann-Whitney $U$ statistic) recovers its well-known consistency as
long as $\mathbb P(X<Y)\neq 1/2$, and remarkable Pitman efficiencies compared to the
$T$-test under location families $g_n(x) = f(x - \theta/\sqrt{n})$. For scale
families $g_n(x) = f(x/\theta_n)/\theta_n$ however, performance relative to the
relevant $F$-test is modest (indeed, many choices of $f$ render $\| S_{n,k} \|_{1,
(0,...,k)/k}$ inconsistent), naturally raising the question of whether a suitable
choice of $w$ may transfer much of $\| S_{n,k} \|_{1,(0,...,k)/k}$'s power in the 
location setting to that of scale families.

\begin{example}[Detecting heteroskedasticity] \label{example:heteroskedasticity}

    \begin{figure}[t] \centering
       \includegraphics[width=\textwidth]{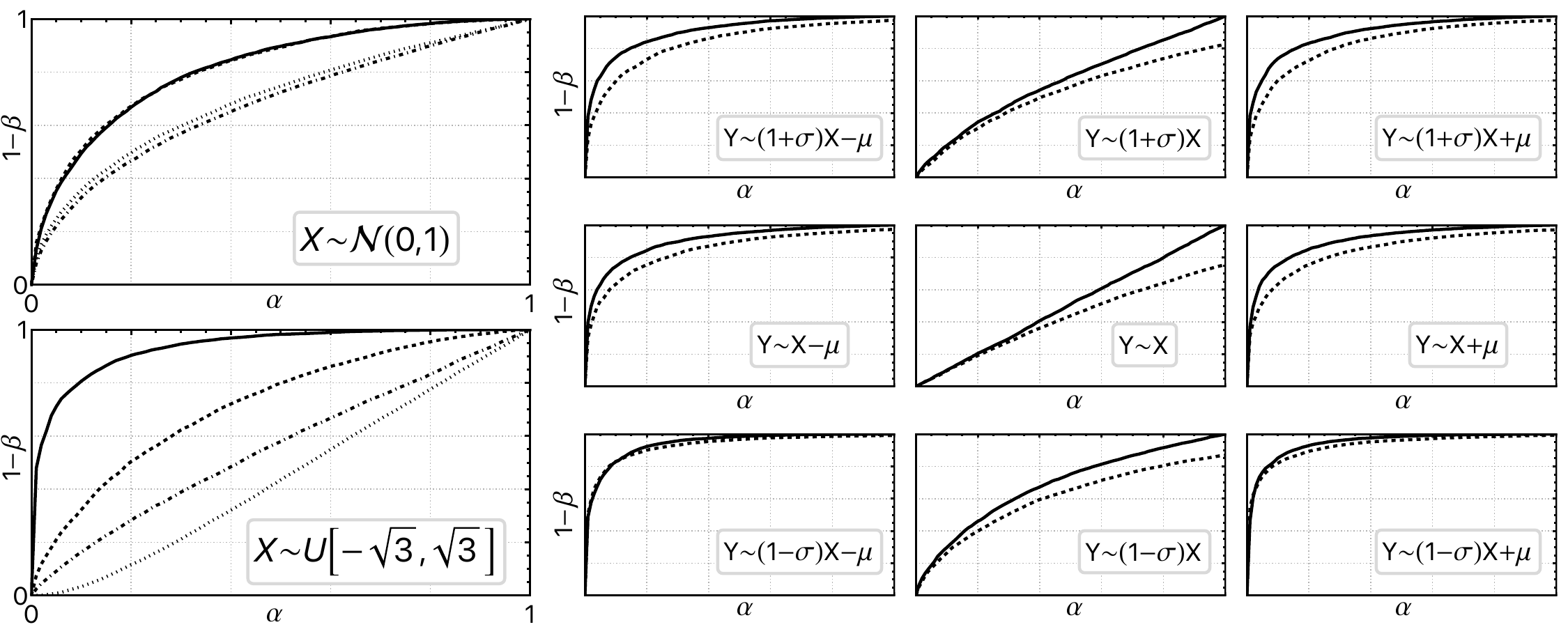}
       \caption{
           \textbf{Left}: ROC curves comparing asymptotic power of $\| S_{n,k}
           \|_{1,w}^1$ as described in Example \ref{example:heteroskedasticity} (solid
           line), the $F$-test (dashed), Bartlett's test (dotted), and the
           Brown-Forsythe test (dot-dashed) in the context of scale shifts
           from $F = \mathcal N(0,1)$ and $F = \operatorname{Uniform\left(
           \left[-\sqrt{3}, \sqrt{3} \right] \right)}$. Simulations used $k=10^5,
           n/k=3/4$ and $\theta_n = \sqrt{1+n^{-1/2}}$ at $10^5$ Monte-Carlo
           iterations. \textbf{Right}: ROC curves comparing power of combined tests
           $\left( \| S_{n,k} \|_{1,w_1}^1, \| S_{n,k} \|_{1,w_2}^1 \right)$
           (solid) and $\left( F\text{-test}, t\text{-test} \right)$ (dashed) as
           described in Example \ref{example:general}. For all panels, $X\sim
           \mathrm{Laplace}\left( 0, 2^{-1/2} \right), n = 18, k = 24,
           \mu=\sigma=3/\sqrt{n}$.
       }
       \label{fig:joint_power} 
    \end{figure}

    Repeating the above calculations for the choice $w_j = w(j/k)$ where 
    \[
        \sqrt{2} w(x) = -\frac{1}{2} + \frac{d}{dt} \lim_{n\to\infty} \sqrt{n}\left[
        \Phi\left( \frac{\Phi^{-1}(x)}{\sqrt{1+n^{-1/2}}} \right) - t\right] =
        \Phi^{-1}(x)^2-1,
    \]
    with $\Phi$ and $\Phi^-1$ the \textnormal{CDF} and inverse \textnormal{CDF} of a
    $\mathcal N(0,1)$ variable respectively, shows that the so obtained $\| S_{n,k}
    \|_{1,w}^1$ is consistent against $G$ as long as $\int_0^1 w(x) g^F(x) \
    \mathrm{d}x \neq 1/2$, which by symmetry of $w$ is the case under, e.g., the
    above-mentioned family of scale shifts if $F$ is symmetric, and provides the
    suitable spacing analogue of Klotz' Gaussian score statistic
    \autocite{klotz1962nonparametric}. $\mu_h$ can be explicitly computed for various
    choices of $F$, and while the corresponding Pitman efficiencies do not mirror 
    uniform lower bounds as were present in the context of location families (indeed,
    \autocite{klotz1962nonparametric} already showed that efficiencies can drop as low as
    $0.47$, conjecturing that arbitrarily small values are possible; which 
    Proposition \ref{prop:klotz} above confirms), they behave favorably for many $F$
    commonly encountered in practice. Some such efficiencies for various choices
    of $F$ are displayed in the left half of Figure \ref{fig:joint_power} (which also
    includes ROC curves for two popular alternatives to the $F$-test: Bartlett's test
    and the Brown-Forsythe test) and Table \ref{tab:pitmans}. The results broadly
    mirror the corresponding efficiencies for Mann-Whitney's $U$ against the $T$-test
    under location shifts, and together with the fact that both the distribution
    exhibited in \autocite{klotz1962nonparametric}  as well as Proposition
    \ref{prop:klotz} leading to efficiencies $<1$ are comparatively impractical,
    render $\| S_{n,k} \|_{1,w}^1$ with this choice of
    $w$ a promising, simple to use candidate for testing against scale shifts or
    variance differences more generally, when Gaussian or other parametric
    assumptions are not available (in order to calibrate the $F$-test, even
    asymptotically, the variance and central fourth moment of $F$ need to be known).
    \begin{table}
        \centering
        \begin{tabular}{l | c c c c c c}
            \toprule
            $F$ & $\quad$ Gaussian & Laplace & Student's $t$ ($\nu = 5$) & Gumbel &
            Cauchy &
            $f_Y \in \mathcal C([a,b]) \setminus \mathcal C_0([a,b])$
            \\\midrule 
            ARE $\,$ & $\quad 1$ & $\approx 1.23$
            & $\approx 2.32$ & $2.59$ & $\infty$ & $\infty$ \\\bottomrule 
    \end{tabular}
        \caption{Asymptotic Relative (Pitman) Efficiency of $\| S_{n,k} \|_{1,w}^2$ as
        described in Example \ref{example:heteroskedasticity} compared to $F$-test under
        various choices of $F$ (the distribution of $X$). The $f_Y \in \mathcal C([a,b])
        \setminus \mathcal C_0([a,b])$ column includes random variables $Y$ whose
        densities take non-zero values at at least one boundary of the support of $Y$.
        Either $a$ and $b$ (not necessarily both) can be finite.
        }
        \label{tab:pitmans}
    \end{table}
\end{example}

\begin{example}[Detecting location and scale shifts]\label{example:general}
    Given the complementary nature of Mann-Whitney's $U$ and the spacing statistic
    discussed in Example \ref{example:heteroskedasticity}, it is appealing to combine
    both tests along the lines of the discussion surrounding Proposition
    \ref{prop:multivariate_extension} of the main
    article. This is possible, since both Mann-Whitney's $U$ and the $\| S_{n.k}
    \|_{1,w}$ statistic discussed in Example \ref{example:heteroskedasticity} belong
    to the same family of spacing statistics, allowing efficient computation and
    inversion of the Laplace transform and joint moments. Such combination may serve
    as a fully non-parametric yet powerful alternative to, e.g., the common practice
    of performing both $F$- and $T$-tests sequentially and reporting Bonferroni- or
    otherwise corrected $p$-values. Asymptotically as $n,k\to\infty$ at similar rate,
    $(\| S_{n,k}/n \|_{1,w_1},\|
    S_{n,k}/n \|_{1,w_2})$ (suitably normalized) is jointly Gaussian with covariance 
    $\int w_1w_2$, and so explicit joint Laplace transforms or
    moment computations are not necessary. However, if sample sizes are small,
    accounting for correlations through such explicit computation may be expected to
    improve power over general (often conservative) $p$-value correction schemes in
    addition to the marginal increases of power between $\| S_{n,k} \|_{1,
    (\Phi^{-1}(j/k)^2-1)/\sqrt{2}}$ and the $F$-test, and $\| S_{n,k}
    \|_{1,(0,...,k)/k}$ and the $T$-test. We illustrate this on the case of comparing
    $( \| S_{n,k} \|_{1,(0,...,k)/k}, \| S_{n,k} \|_{1,(\Phi^{-1}(j/k)^2-1)/\sqrt{2}}
    )$ against the $( \text{F-Test}, \text{T-Test})$ combination in the right half of
    Figure \ref{fig:joint_power}, where ROC curves for mixtures of location- and
    scale-shifts are plotted anchored on $X\sim \mathrm{Laplace}\left( 0, 1/\sqrt{2}
    \right)$.
\end{example}

\begin{figure}[t] \centering
    \includegraphics[width=\textwidth]{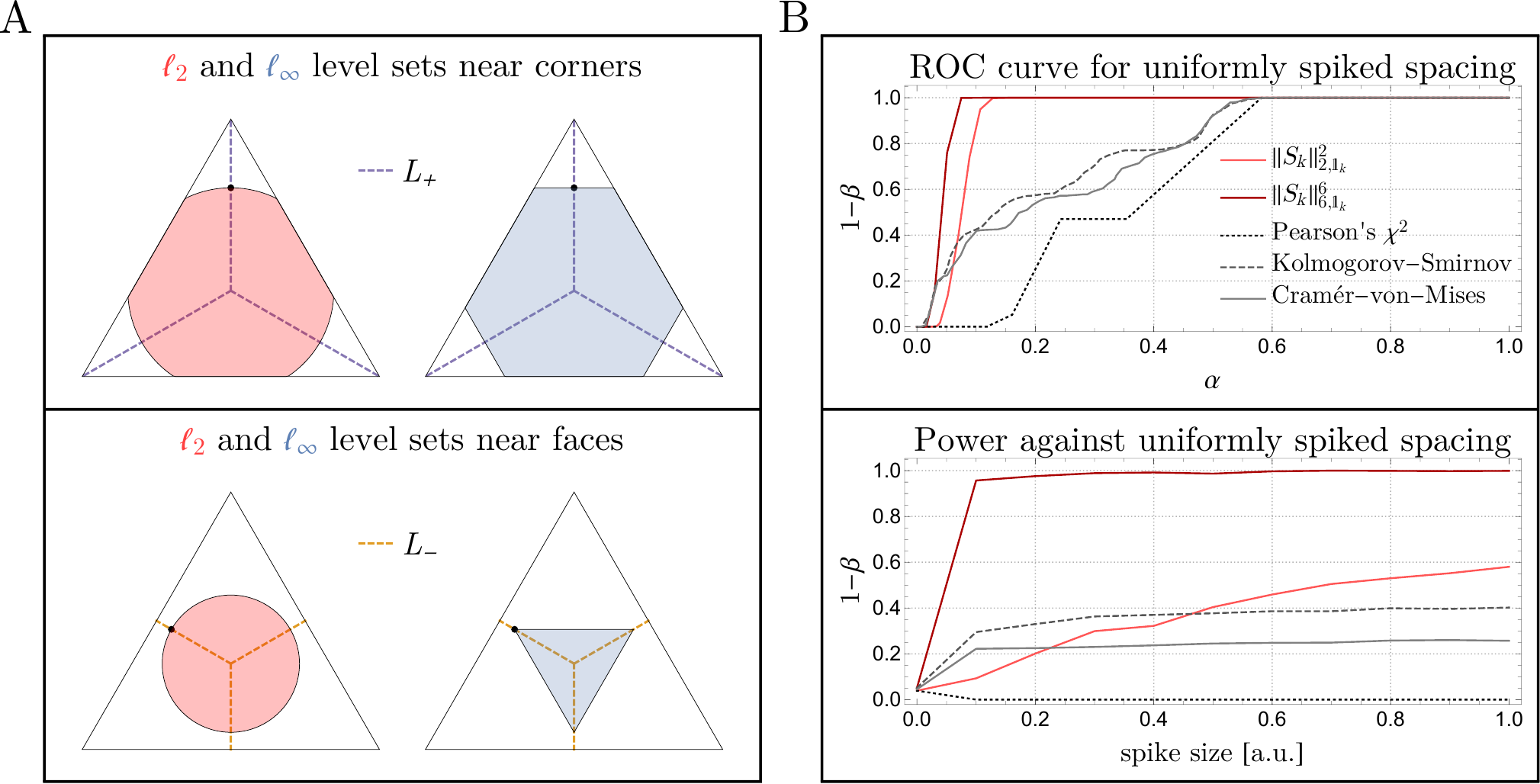}
    \caption{
        Analysis of spiked spacing model (described in Example~\ref{example:spiked}). \textbf{A.}
        Illustration of tail probabilities on $\Delta^{2}$ in the cases of $p = 2$ and $p = \infty$,
        and the samples (denoted by solid dots) giving rise to them.  While a fixed sample near line
        segments in $L_+$ (purple, dashed lines in top panel) produces smaller sub-level sets in
        $\ell_2$ than $\ell_{\infty}$ (thereby increasing the $p$-value of said sample, which
        corresponds to $1$ minus the shaded area), this trend  reverses for observations near
        line segments in $L_{-}$ (orange, dashed lines in bottom panel).
        \textbf{B.} ROC and power curves. The spiked spacing model largely
        concentrates around $L_{\text{corner}}$ in $\Delta^{k-1}$, with the degree of this
        concentration increasing with spike size. As a consequence, $p$-norms of samples generated
        under such alternative tend to separate more markedly for larger $p$, which in turn affords
        increases in power of $\|\bfS_{k}\|_{p,\bfmath{1}_k}^p$ when $p>2$. The experimental design
        and choices of under- and overdispersed distributions follows that of Figure $1$ in the main
        article.
    } 
    \label{fig:spiked_spacing}
\end{figure}

\begin{example}[Spiked spacing model and higher $p$-norms]\label{example:spiked}
    Due to its correspondence with the likelihood-ratio test, the choice of $p=1$ is optimal
    whenever data arrives in an $iid$ fashion as was the case in Examples
    \ref{example:heteroskedasticity} and \ref{example:general}. When observations exhibit
    correlation or are otherwise structured, larger values of $p$ may become relevant. We illustrate
    this phenomenon by revisiting the one-sample test in (Ref),
    where we sought to distinguish exponential arrival times from over- or
    underdispersed alternatives. We consider an alternative hypothesis $G$ of the joint
    distribution of $T_1, \dots, T_k$ that is both  over- and underdispersed
    in the following sense: Under $G$, arrival times are again drawn $iid$ from an
    underdispersed distribution $G_{1}$, with the exception of a single randomly chosen $T_K$ (i.e., $K \sim
    \mathrm{Uniform}([k])$) whose law $G_{2}$ now exhibits overdispersion. We call this
    overdispersed $T_K$ the spiked or outlier arrival time. Though the subsequent analysis is
    phrased in terms of this spiked spacing model, much of its reasoning pertains to similar outlier
    or correlation models of this kind as well. 

    To design a test that reliably detects this spiked spacing model, we first observe that the
    symmetry in $T_1, \dots, T_k$ suggests little benefit of choices for $\bfw$ other than
    $\bfmath{1}_k$, leaving $p$ as the sole parameter to optimize. It is clear that on the level of
    normalized spacings, the null and alternative distributions differ only by the presence of
    exactly one particularly large segment, the index of which is random, and so a generalized
    likelihood ratio test is effectively based on the length of the longest segment. In terms of $\|
    \bfS_{n,k} \|_{p,w}^p$ this corresponds to the choice $p=\infty$. To reason about intermediate
    values of $p$ between $1$ and $\infty$, it is useful to clarify and compare the geometry that
    various $\ell_p$ balls give rise to when intersected with $\Delta^{k-1}$: as the $2$-dimensional
    illustrations of Figure~\ref{fig:spiked_spacing}A demonstrate, the (normalized)
    intersection volume $V_k^p(\bfs) = \mu_{\Delta^{k-1}} ( \|\bfS_k\|_{p,\bfmath{1}_k}^p \leq
    \|\bfs\|_{p,\bfmath{1}_k}^p )$ depends on the precise location of the observation
    $\bfs$. If $\bfs$ lies exactly on any of the line segments $L_{\text{corner}} =
    \Big\{\lineseg{\frac{1}{k}\bfmath{1}_k}{\bfe_i}\Big\}_{i\in[k]}$, where $\bfe_i$ is the
    $i^{\text{th}}$ standard basis vector, then $V_k^p(\bfs) \subset V_k^q(\bfs)$ whenever $p < q$,
    while $V_k^p(\bfs) \supset V_k^q(\bfs)$ in case $\bfs$ falls precisely on any of the line
    segments $L_{\text{mid}}=\Big\{ \lineseg{\frac{1}{k}\bfmath{1}_k}{\bfm_i} \Big\}_{i\in[k]}$, where
    $\bfm_i = \frac{1}{k-1}(\bfmath{1}_k - \bfe_i)$ is the midpoint of the $(k-2)$-dimensional face
    opposite of vertex $\bfe_i$. Since $p$-values are $1-V_k^p(\bfs)$, it follows that
    tests based on $\|\bfS_k\|_{\infty,\bfmath{1}_k}^{\infty}$ should be most powerful in the former
    scenario, while $\|\bfS_k\|_{2,\bfmath{1}_k}^2$-based tests benefit from the latter scenario, with
    intermediate localizations giving rise to optimal $p^*$ between $2$ and $\infty$. In the spiked
    spacing model, the support of $G$ centers around the line segments $L_{\text{corner}}$, and so
    we expect choices of $p$ larger than $2$ to be profitable. Indeed, carrying out simulations as
    in Figure~\ref{fig:spiked_spacing}B reveals this to be true, with precise values of $p^*$
    depending on the distributional details $G_{1}$ and $G_{2}$.  Generally, $p^*$ is attained
    around $4$ or $5$ for modest amplitudes of the spiked $T_K$ and/or moderate degrees of
    underdispersion in the remaining arrival times, and stabilizes at $6$ for more pronounced levels
    of spiking and/or underdispersion. Past $p = 6$, ROC and power curves tend to change only
    slightly.  
\end{example}

\begin{figure}[t]
    \centering
    \includegraphics[width=\textwidth]{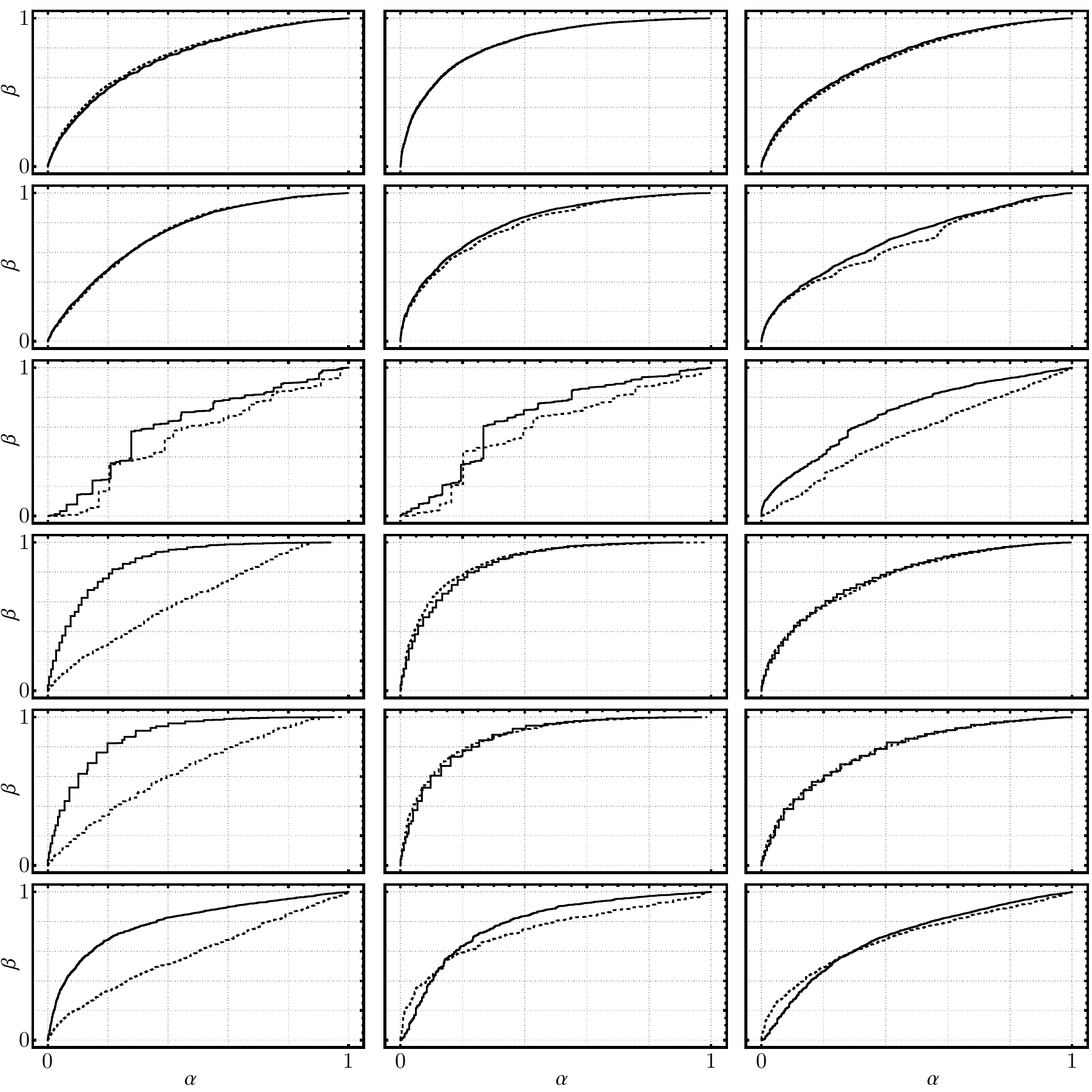}
    \caption{Comparison of ROC curves corresponding to
    $\|S_{n,k}\|_{1,w}$ (solid lines) and $T_{n,k}$ (dashed) for $n=10,
    k=5$, and various choices of $w, F$, and $G$. The three columns are
    associated with null distributions following
    $\operatorname{Uniform}[0,1], \mathcal N(0,1)$ and
    $\operatorname{Cauchy}(0,1)$, while the first and second set of
    three rows correspond to $Y \sim X+1$, and $Y\sim 1.1\cdot X$,
    respectively. Choices of $w$ are fixed row-wise, and read $w_1 =
    \Phi^{-1}, w_2(x) = x^6, w_3 = r_{\mu}, w_4 =
    \left(\Phi^{-1}\right)^2, w_5(x) = (x-1/2)^2, w_6 = r_{\sigma}$,
    where $r_{\mu}$ and $r_{\sigma}$ are densities obtained from
    normalizing a Brownian motion on $[0,1]$.
    }
    \label{fig:roc}
\end{figure}

\newpage

\begin{figure}[t]
    \centering
    \includegraphics[width=\textwidth]{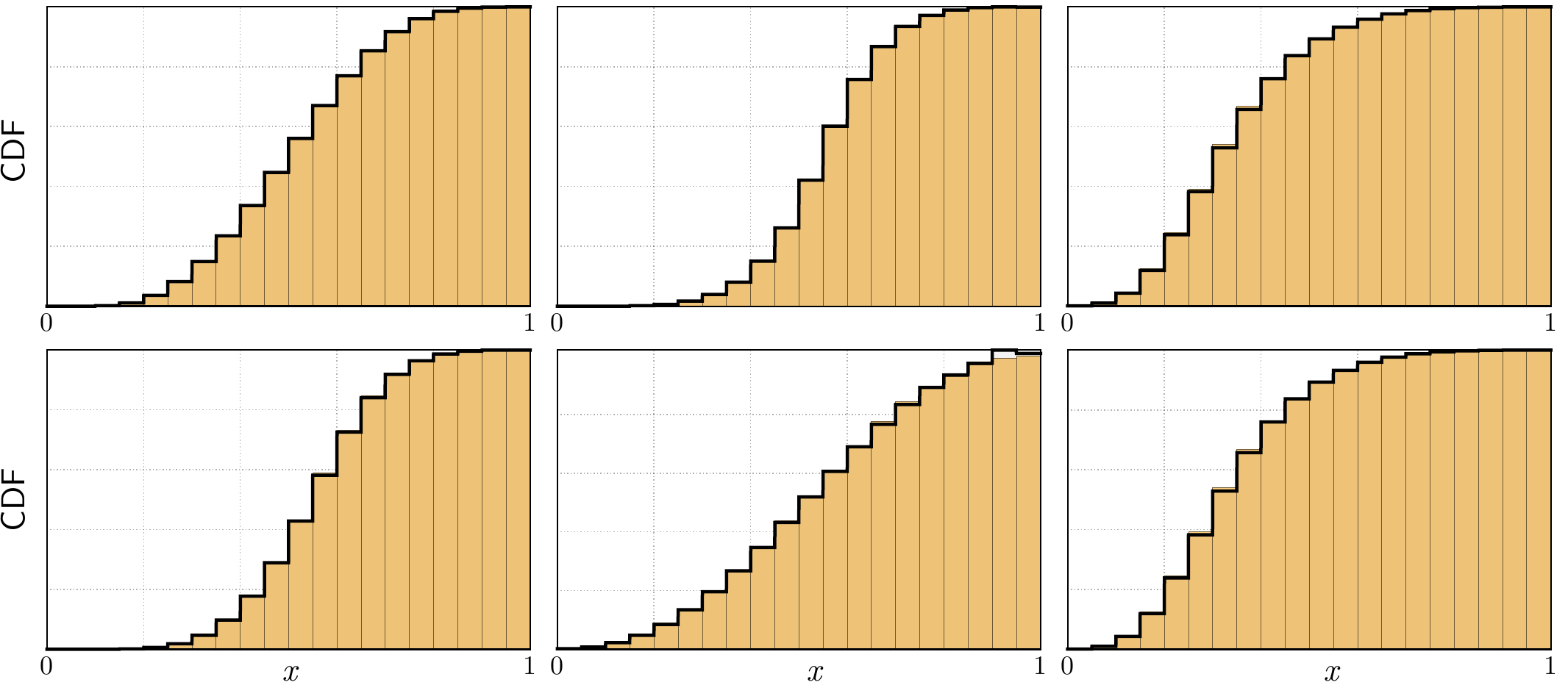}
    \caption{Comparison of CDFs of $\| S_{n,k} \|_{1,w}$ obtained from numerical inversion (solid line) of its Laplace transform (cf. Theorem 1 in the main article) with simulations (bar chart). Six choices of $w \in \mathbb R^7$ were sampled by drawing uniform $[0,1]$ entries i.i.d. and normalizing to $\| w \|_{\infty} = 1$, and corresponding draws of $\| S_{6,15} \|_{1,w}$ simulated $10,000$ times. 
    }
    \label{fig:cdf}
\end{figure}

\end{document}